\documentclass[12pt]{article}
\usepackage[psamsfonts]{amssymb}
\usepackage{amsfonts}
\usepackage[mathcal]{eucal}
\newcommand{\g}{\mathfrak{g}}
\renewcommand{\a}{\mathfrak{a}}
\renewcommand{\b}{\mathfrak{b}}
\newcommand{\h}{\mathfrak{h}}
\renewcommand{\k}{\mathfrak{k}}
\newcommand{\q}{\mathfrak{q}}
\newcommand{\p}{\mathfrak{p}}
\newcommand{\n}{\mathfrak{n}}
\newcommand{\m}{\mathfrak{m}}
\renewcommand{\l}{\mathfrak{l}}
\renewcommand{\j}{\mathfrak{j}}
\newcommand{\s}{\mathfrak{s}}
\renewcommand{\t}{\mathfrak{t}}
\newcommand{\z}{\mathfrak{z}}

\newcommand{\Q}{\mathbb{Q}}
\newcommand{\Z}{\mathbb{Z}}
\newcommand{\R}{\mathbb{R}}
\newcommand{\C}{\mathbb{C}}

\newcommand{\qed}{{\null\hfill\ \raise3pt\hbox{\framebox[0.1in]{}}}\break\null}
\newtheorem{theo}{Th\'eor\`eme}
\newtheorem{prop}{Proposition}
\newtheorem{lem}{Lemme}
\newtheorem{cor}{Corollaire}
\newtheorem{rem}{Remarque}
\newtheorem{defi}{D\'efinition}
\newcommand{\ste}{\hfill\break}

\begin{document}
\def\b{\mathfrak{b}}
\def\c{\mathfrak{c}}
\def\d{\mathfrak{d}}
\def\e{\mathfrak{e}}
\def\f{\mathfrak{f}}
\def\h{\mathfrak{h}}
\def\i{\mathfrak{i}}
\def\j{\mathfrak{j}}
\def\r{\mathfrak{r}}
\def\k{\mathfrak{k}}
\def\q{\mathfrak{q}}
\def\p{\mathfrak{p}}
\def\n{\mathfrak{n}}
\def\m{\mathfrak{m}}
\def\l{\mathfrak{l}}
\def\j{\mathfrak{j}}
\def\s{\mathfrak{s}}
\def\t{\mathfrak{t}}

\def\u{\mathfrak{u}}
\def\v{\mathfrak{v}}
\def\w{\mathfrak{w}}
\def\x{\mathfrak{x}}
\def\K{\mathbb{K}}

\def\dem{ {\em D\'emonstration : \ste }}
\def\sac{ sous-alg\`ebre de Cartan }
\def\beq{\begin{equation}}
\def\eeq{\end{equation}}
\def\z{\mathfrak{z}}
\textheight 659pt\textwidth  444pt
\oddsidemargin  -1mm
\evensidemargin -1mm
\topmargin      -8mm
\pagestyle{plain}

\title{\bf Sur les triples de Manin pour les alg\`ebres de Lie
r\'eductives complexes
}
\author{Patrick Delorme}

\maketitle
\setcounter{section}{-1}

{\bf Abstract}\ste
We study Manin triples for a reductive Lie algebra, $\g$. First, we
generalize results of  E. Karolinsky, on the classification of
Lagrangian subalgebras (cf. KAROLINSKY E., {\em A Classification of
Poisson homogeneous spaces of a compact Poisson Lie group},
Dokl. Ak. Nauk, 359 (1998), 13-15). Then we show that, if $\g$ is
non commutative, one can attach , to each Manin triple in
$\g$ , another one for a  strictly smaller reductive
complex Lie subalgebra of $\g$. We study also the inverse process.
\ste
\section{Introduction}

Let $\g$ be a finite dimensional, complex,  reductive Lie algebra.
One says that a  symmetric, $\g$-invariant,
 $\R$-bilinear form on $\g$ is a {\bf Manin form} if and only if
its signature is
$(dim_{\C}\g, dim_{\C}\g)$. A Manin form is called special, if its
restriction to any complex semisimple Lie algebra is non zero. It is
easily seen that for a complex semisimple (resp. simple) Lie algebra,
every non degenerate (resp. non  zero)  symmetric,
 $\g$ invariant,
 $\R$-bilinear form on $\g$ is a Manin form (resp. special). \ste
Recall that a {\bf Manin-triple} in $\g$ is a triple $(B,\i, \i')$,
where
$B$ is a Manin form and  where
$\i,
\i'$ are real Lie subalgebras of $\g$, isotropic for $B$, and such
that $\i+ \i'=\g$. Then this is a direct sum and $\i, \i'$ are of
real dimension equal to the complex dimension of $\g$.
\ste    Our goal is to construct all Manin-triples of $\g$, up to
conjugacy under the action on
$\g$ of the connected, simply connected Lie group $G$ with Lie
algebra $\g$, by induction on the rank of the derived algebra of
$\g$.\ste
One calls {\bf af-involution } of a complex semi-simple Lie algebra
$\m$, any involutive automorphism of $\m$, $\sigma$ such that
 there exists :\ste 1) an ideal ${\tilde  \m_{0}}$ of $\m$, and a
real form
${\tilde \h}_{0}$ of ${\tilde \m_{0}},$\ste   2) simple  ideals of
$\m$, $\m'_{j}$,
$\m''_{j}$, $j=1,\dots,p$,
\ste
3)  an  isomorphism of   complex Lie algebras, $\tau_{j}$,
between
$\m'_{j}$ and
$\m''_{j}$, $j=1,\dots,p$, \ste
such that, denoting by  $ \h_{j}:=\{(X, \tau_{j}(X)\vert X
\in
\m'_{j}\}$, and by  $\h$  the fixed point set  of $\sigma$,
one has  :

$$\m={\tilde \m_{0}}\oplus (\oplus_{ j=1,\dots,p} (\m'_{j}\oplus
\m''_{j}))$$
$$\h= {\tilde \h}_{0}\oplus (\oplus_{ j=1,\dots,p} \h_{j}) $$
Notice that an $\R$-linear involutive automorphism of $\m$ is
determined by its fixed point set, as the set of antiinvariant
points is just the orhogonal of the fixed point set,  for the
Killing form of
$\m$, viewed as a real Lie algebra.
The following Theorem generalize previous results of E.Karolinsky
(cf [K2], Theorem 1 (i) and [K1] for the proof, see also [K3]
Proposition 3.1), as we do not make any restriction on the Manin
form.
 \ste{\bf Theorem 1}
\ste{\em Let
$B$ be a Manin form and let
$\i$ be a {\bf Lagrangian subalgebra of $\g$ for $B$}, i.e. a real
Lie subalgebra of $\g$ of dimension the complex dimension of $\g$
and isotropic for $B$ . Then :
\ste a) If we denote by $\p$ the normalizer in $\g$ of the nilpotent
radical,
$\n$, of $\i$, then $\p$ is a parabolic subalgebra of $\g$, with
nilpotent radical equals to $\n$.\ste
b) Let $\l$ be a Levi subalgebra of $\p$ (i.e. $\l$ is a reductive
Lie subalgebra of $\p$ with $\p=\l\oplus\n$), and denote  by $\m$
its
derived ideal and $\a$ its center,then  the intersection, $\h$, of
$\i$ and $ \m$ is the fixed point set of an af-involution of
$\m$, which is  isotropic for
$B$. Moreover, if  $B$ is special $\h$ is  a real form of $\m$\ste
c) The intersection,
$\i_{\a}$, of
$\a$ and
$\i$, is isotropic for
$B$, and its real dimension equals  the complex dimension of
$\a$.\ste  d) One has
$\i=\h\oplus\i_{\a}\oplus\n$.\ste  Reciprocally, any real   Lie
subalgebra, $\i$, of $\g$, which is of this form is
Lagrangian  for
$B$. \ste Then,  one
says that {\bf
$\i$ is under $\p$ }}
\ste
One chooses a Cartan subalgebra $\j_{0}$ of $\g$, and a Borel
subalgebra of $\g$ containing $\j_{0}$, $\b_{0}$. Then, from
Theorem 1 and the Bruhat decomposition, one sees (cf. Proposition
1) that  every Manin triple  is conjugated, under
$G$, to a Manin
triple
$(B,\i, \i')$ such that  $\i$ is under $\p$ and  $\i'$ is under
$\p'$, with $\p$ containing $\b_{0}$ and $\p'$ containing the
opposite Borel subalgebra to
$\b_{0}$, with respect to $\j_{0}$. A Manin
triple  satisfying these
conditions will be called {\bf standard}, under $(\p, \p')$.
\ste
If $\r$ is a real subalgebra of $\g$, we denote by $R$ the
analytic subgroup of
$G$ with
 with Lie algebra $\r$ .\ste
If $\e$ is an abelian real subalgebra of $\g$, $\r$ is  an
$\e$-invariant subspace of
$\g$, let
$\Delta(\r,\e)$ the set of weights of $\e$ in $\r$, which is the
subset of $Hom_{\R}(\e,\C)$ whose elements are joint eigenvalues of
elements in $\e$ acting on $\r$.The weight space of
$\alpha$ in
$\Delta(\r,\e)$ is denoted by $\r^{\alpha}$. \ste
Let $\p$, $\p'$ be
given as above,  and let
$B$ be a Manin form on $\g$.\ste \ste
{\bf Theorem 2}  {\em \ste  If there exists a standard
Manin
triple
$(B,\i,
\i')$  under $(\p, \p')$,
 then $\p$ or
$\p'$ is different from
$\g$ }.\ste\ste {\bf Theorem 3} {\em \ste Suppose that $(B,\i,
\i')$ is a standard  Manin triple under $(\p,
\p')$). Let  $p^{\n'}$  be the projection of $\g$ on the
$\j_0$-invariant supplementary subspace of the nilpotent radical
$\n'$ of $\p'$ , with kernel $\n'$ . Let $\l\oplus \n$  the
Langlands decomposition of   $\p$, such that $\l$ contains
$\j_{0}$. Set
$\i_{1}=p^{\n'}(\tilde{\h}\cap \p')$, where $\tilde{\h}=\i\cap \l$.
One defines similarly $\i'_{1}$. \ste
Then $\i_{1},\i'_{1}$ are contained in $\l\cap \l'$. Moreover, if
$B_{1}$ denotes the restriction of $B$ to $\l\cap \l'$,
$(B_{1},\i_{1},\i'_{1})$ is a Manin triple for $ \l\cap \l'$. We
set $\g_{1}:=\l\cap\l '$.\ste We will use freely the notation of
Theorem 1 for $(B_{1},\i_{1},\i'_{1})$, which is called
the {\bf predecessor  of the standard Manin triple }$(B,\i,
\i')$ .\ste} \ste
{\bf Theorem 4} {\em \ste
(i) Every Manin triple under $(\p,\p')$ is conjugate to
a standard Manin triple as above, $(B,\i, \i')$, with predecessor
$(B_{1},\i_{1},\i'_{1})$, for which there exists a fundamental
Cartan subalgebra $\tilde{\f}$ (resp.
$\tilde{\f'}$) of
$\i$ (resp.$\i'$),
included in $\i_1$ (resp. $\i'_1$), such that, if we denote by
$\sigma$ (resp. $\sigma '$)  the af-involution of $\m$
(resp. $\m'$), with fixed point set $\h$ (resp. $\h'$),
then the pair $(\tilde{\f}, \sigma)$
satisfies the following properties (the pair $(\tilde{\f'},
\sigma')$ satisfying similar properties with the obvious
changes) :\ste 1)
The map $\sigma$ is an af-involution of
$\m$ with fixed point set $\h$, and $\tilde{\f}$ is a fundamental
Cartan subalgebra of
$\i_{1}$ such that $\tilde{\f}= \f+\i_{\a}$, where
$\f=\tilde{\f}\cap \h$ is a fundamental Cartan subalgebra of $\h$
and
$\i_{\a}=\i
\cap
\a$. Moreover
$\i_{\a}$ is Lagrangian  for $B$ restricted to  $\a$, and $\h$ is
isotropic for
$B$.
Let us denote by $\j$ the centralizer of ${\tilde \f}$ in $\g$.
This is a Cartan subalgebra of $\g$, contained in $\l\cap\l'$. \ste
2) The intersection of
$\h$ with
$\n'$ is reduced to zero.
\ste 3)There  exists a Borel subalgebra, $\b$, of  $\m$,
containing  $\j\cap\m$ and included in  $\m\cap \p'$, such that
$\sigma(\b)+\b=\m$.\ste
4) There exists a unique  Langlands
decomposition of
$\p_1$ (where
$\i_1$ is under $\p_1$),  ${\underline
\l_{1}}\oplus \n_{1}$, with ${\underline \l_{1}}$ containing
${\tilde
\f}$. Then
${\underline \m_{1}}:=[{\underline \l_{1}},{\underline \l_{1}}]$
is equal to the derived ideal of
$(\m\cap\l')\cap
\sigma(\m\cap\l')$.\ste
5) If $\alpha\in \Delta (\m, \j)$, one can define
${\underline \alpha}\in
\Delta (\m, \j)$ by the condition : $\sigma(\m^{\alpha})=
\m^{{\underline \alpha}}$.\ste One has :
$\n_{1}=\oplus_{\alpha \in \Delta(\n'\cap \l,
\j)\cap\underline{\Delta(\l\cap\l',
\j)}}\m^{{\underline \alpha}}$.\ste 6) The restriction of
$\sigma$ to
${\underline \m_{1}}$ is equal to the af-involution of ${\underline
\m_{1}}$ having
${\underline \h_{1}}:=\i_{1}\cap{\underline \m_{1}}$ as fixed point
set.
\ste
 One says
that the standard  Manin triple
$(B,\i,\i')$ is {\bf linked }, or linked to its  predecessor
$(B_{1},\i_{1},\i'_{1})$  ,
with  {\bf link }
$(\tilde{\f},\tilde{\f'})$.\ste
(ii) Reciprocally, if $B$ is a Manin  form, if $(B_{1}, \i_{1},
\i'_{1})$ is a
  Manin
triple  for $\l\cap\l'$, and if $(\tilde{\f},\sigma)$
(resp.
$(\tilde{\f'},\sigma ')$) satisfies the conditions 1) to 6) (resp.
the conditions 1) to 6) with the obvious changes), then the triple
$(B,\underline{\i},\underline{\i'})$, where
$\underline{\i}:=\h+\i_{\a}+\n$, and
$\underline{\i'}$ is defined similarly, is a standard  Manin
triple , under $(\p,\p')$, linked to
$(B_{1},\i_{1},\i'_{1})$
with link
$(\tilde{\f},\tilde{\f'})$.}
\ste\ste  One shows easily that every Manin triple for $\g$ is
conjugated under $G$ to a standard linked Manin triple  having a
sequence of predecessors which are standard and linked (cf.
Proposition 5). We say that such a Manin
triple  is {\bf strongly standard}. Then the sequence is
stationnary and its ending element,denoted by $(B_{p},
\i_{p},\i'_{p})$ (where $p$ denotes the first index
 for which $\g_{p}=\j_{0}$) is a Manin triple in
$\j_{0}$, where
$B_{p}$ is the restriction of $B$ to
$\j_{0}$. It depends only on the conjugacy class under $G$ of the
original Manin triple and it is called the  {\bf  bottom of this
conjugacy class}.  Then, if  $(B_{1},\i_{1},\i'_{1})$ is as in
Theorem 3 (ii), and moreover strongly standard, the existence of a
fundamental  Cartan subalgebra  $\tilde{\f}$ of
$\i_{1}$ with $\tilde{\f}$ is equal to the sum of its intersections with
$\m$ and $\a$ (as in condition 1) of the Theorem 4) is
equivalent to the fact that $\i_{p}$ satisfies the same property.
Moreover, in this case $\tilde{\f}$ and $\i_{p}$ are conjugate by
an element $x_{1}$ of  $I_{1}I_{2}...I_{p}$. More precisely
$x_{1}=y_{1}\dots y_{p}$, where the $y_{j}$ are elements of
$I_{j}$ such that for all $j=2,\dots ,p$, one has
$Ad(y_{j}y_{j+1}\dots y_{p})(\i_{p})\subset \i_{j}\cap\i_{j-1}$
\ste This allows us to partly reduce the construction of  an
af-involution, for which the fixed point set has
$\f=\tilde{\f}\cap \m$ as fundamental Cartan subalgebra, and
satisfying  the properties 1) to 5),  to a similar
problem with
$ \f$ replaced by
$\i_{p}\cap \m$, ${\underline \l_{1}}$, by $\l_{1}$, where
$\l_{1}$ contains $\j_{0}$, and $\l_{1}\oplus \n_{1}$ is a
Langlands decomposition of $\p$ (cf. Lemme 13). The
correspondance between the solutions of the two problems is given
by conjugacy under $x_{1}$. In principle, this problem can be
solved explicitely , the possibilities for $\i_{p}\cap\m$ being
finite . Up to conjugacy by
$J_{0}$, there are only finetely many  solutions (Propositions 6,
7). The condition 6) adds an extra
restriction.\ste  One wants also to decide when the  Manin triples
that are built with the help of Theorem 4 {\em(ii)} are conjugate
under
$G$.
\ste This question reduces, by induction, to
determining when strongly standard Manin triples,  with a given
predecessor
$ (B_{1},,\i_{1},\i'_{1})$, and which are under $(\p,\p')$ are
conjugate.\ste Let
$\overline{I}$ be the normalizer of $\i$ in $G$. It  is equal to
$\overline{H}AN$, where
$\overline{H}=\{m\in M \vert Ad \>m_{\vert\m}\>\> commutes
\>\>with\>\> \sigma\}$. Then, the problem of conjugacy ultimately
reduces to deciding whether the intersection of $x \overline{I}$
with
$x'\overline{I'}$ is non empty. Here, the   elements $x, x'$ of $G$
can be determined from the data , in particular from the
$x_{1}$ above, corresponding to each of the triples.
\ste {\bf Comment on proofs:}\ste  The proof of Theorem 1  rests on
the properties  of the radical and of the nilpotent radical of a
Lie
algebra (cf. [Bou], Groupes et alg\`ebres de Lie, Chapitre I) , and
the decomposition of elements in a semisimple  Lie algebra as a sum
of commuting nilpotent and semisimple elements. As we said already,
our method of proof is probably a mixing of the Karolinsky's method
for proving his  results quoted above. \ste  In the proof of
Theorem 2, we use a result of Gantmacher [G] which ensures that
every  automorphism of a real  semisimple Lie algebra has  non zero
fixed points. \ste The proof of Theorem 3 is elementary.
\ste In the proof of Theorem 4, we  use  the characterization, by
T. Matsuki ([M1], [M2]) of open orbits in generalized flag manifolds
of a complex semisimple Lie group, under the action of the fixed
point set of an involutive automorphism.

\ste{\bf Aknowledgment} : I thank very much C. Klimcik for
suggesting me this work, and for many interesting discussions. I
thank also J.L. Brylinski for pointing out to me the work of E.
Karolinsky.

\section{Sous-alg\`ebres de Lie  Lagrangiennes}
Dans tout l'article, alg\`ebre de Lie voudra dire alg\`ebre de Lie de
dimension finie.\ste
 Si $\g$ est une alg\`ebre de Lie on notera souvent  $\g^{der}$ son
id\'eal d\'eriv\'e.\ste  Soit $\a$ est une alg\`ebre de Lie ab\'elienne sur
$\K=\R$ ou $\C$, $V$ un $\a$-module complexe. Pour $\lambda\in
Hom_{\K}(\a,\C)$, on note
$V^{\lambda}:=\{v\in V\vert Xv=\lambda(X) v, \>\> X\in \a\}$, qui
est appel\'e le sous-espace de poids $\lambda$ de $V$. On dit que
$\lambda $ est un poids de $\a$ dans $V$ si $V^{\lambda}$ est non
nul et on note $\Delta (V,\a)$ l'ensemble des poids non nuls de
$\a$ dans $V$.\ste
Si $G$ est un groupe de Lie, on notera $G^{0}$ sa composante neutre.

\begin{lem}\ste(i)Soit $\g$ une alg\`ebre de Lie semi-simple complexe,
$\g_{1},\dots, \g_{n}$  ses id\'eaux simples. Toute forme
$\R$-bilin\'eaire $\g$-invariante sur $\g$ est du type $B_{\lambda}$
(ou
$B_{\lambda}^{\g}$), o\`u $\lambda=(\lambda_{1},\dots,
\lambda_{n})\in\C^{n}$ et :
$$B_{\lambda}(X_{1}+\dots+X_{n},Y_{1}+\dots+Y_{n})=
\sum_{i=1,\dots,n}Im(\lambda_{i}K_{\g_{i}}(X_{i},Y_{i})), $$ si
les $X_{i},Y_{i}$ sont des \'el\'ements de $\g_{i}$. Ici $K_{\g_{i}}$
d\'esigne la forme de Killing de $\g_{i}$.\ste
En particulier, une telle forme est sym\'etrique et les id\'eaux
simples sont deux \`a deux orthogonaux pour une telle forme.\ste
(ii) La forme $B_{\lambda}$ est non d\'eg\'en\'er\'ee si et seulement si
chacun des $\lambda_{i}$ est non nul. Elle est alors de signature
$(dim_{\C}\g, dim_{\C}\g)$.\ste
(iii) La restriction de $B_{\lambda}^{\g}$ \`a une sous-alg\`ebre
complexe et simple, $\s$, de $\g$, est de la forme $B^{\s}_{\mu}$,
o\`u $\mu=\sum_{i=1,\dots,n}q_{i}\lambda_{i}$ et pour chaque $i$,
$q_{i}$ est un nombre rationnel positif. De plus $q_{i}$ est non
nul si et seulement si $\s$ a un crochet non nul avec $\g_{i}$.
\end{lem}{\em D\'emonstration : }Traitons d'abord le cas o\`u $\g$ est
simple, auquel cas $n=1$. La donn\'ee d'une forme
$\R$-bilin\'eaire $\g$-invariante sur $\g$, \'equivaut \`a celle d'une
application $\R$-lin\'eaire entre $\g$ et $Hom_{\R}(\g, \R)$, qui
commute \`a l'action de $\g$, regard\'ee comme alg\`ebre de Lie r\'eelle.
Mais la partie imaginaire de la forme de Killing de $\g$ (regard\'ee
comme complexe) d\'etermine un isomorphisme de $\g$-modules entre
$\g$ et $Hom_{\R}(\g, \R)$. Finalement, la donn\'ee de notre forme
\'equivaut \`a la donn\'ee d'un endomorphisme, $T$, $\R$-lin\'eaire de
$\g$, commutant \`a l'action de $\g$.
\ste Soit $\k$ une forme r\'elle compacte de $\g$. On \'ecrit :$$
T(X)=Re\>T(X)+i Im \>T(X),\>\>X\in \k,$$ o\`u $Re\>T(X),Im\>T(X)\in
\k$. Alors
$Re\>T$, $Im \>T$ sont des \'el\'ements de $Hom_{\k}(\k,\k)$. Comme
$\k$ est simple, il r\'esulte du lemme de Schur que cet espace est
\'egal \`a $\R \> Id_{\k}$. Donc, il existe $\lambda\in \C$ tel que
:$$T(X)=\lambda X,\>\> X\in \k$$ Maintenant, si $X,Y\in \k$, on a :
$$T(i[X,Y])=T([iX,Y])=[iX,TY]$$
la derni\`ere \'egalit\'e r\'esultant du fait que $T$ commute \`a  l'action
de $\g$. Joint \`a ce qui pr\'ec\`ede, cela donne :
$$T(i[X,Y])=\lambda i[X,Y], \>\>X,Y\in \k$$
Comme $[\k,\k]=\k$, on conclut que $T$ est la multiplication par
$\lambda$. D'o\`u (i) dans le cas o\`u $\g$ est simple.\ste Supposons
maintenant que $\g$ soit la somme directe de deux id\'eaux
$\g'$, $\g''$. Soit $B$ une  forme
$\R$-bilin\'eaire $\g$-invariante sur $\g$. On a :
$$B([X',Y'],X'')=-B(Y',[X',X''])=0, \>\>X',\>\> Y'\in\g'
,\>\>X''\in \g'', $$ la derni\`ere \'egalit\'e r\'esultant du fait
que $\g'$,
$\g''$ commutent entre eux. Comme $[\g,\g]=\g$, on a aussi
$[\g',\g']=\g'$. Finalement $\g'$ et $\g''$ sont orthogonaux.
Alors,  on d\'eduit (i) pour $\g$ semi-simple du cas o\`u $\g$ est
simple.\ste (ii) L'assertion sur la non nullit\'e des $\lambda_{i}$
est
claire. Pour l'\'etude de la signature, on se ram\`ene au cas o\`u $\g$
est simple. Supposons $\lambda\in\C$, non nul. Soit $\k$ une forme
r\'eelle compacte de $\g$. On fixe une base $X_{1},\dots, X_{l}$  de
$\k$. On choisit une racine carr\'ee $\mu$ de
$i\lambda^{-1}$.On pose $Y_{i}=\mu X_{i}$, $Z_{i}=i\mu X_{i}$.
Alors
 $B_{\lambda}(Y_{i},Z_{j})=0$,
$B_{\lambda}(Y_{i},Y_{j})=\delta_{i,j}$,
$B_{\lambda}(Z_{i},Z_{j})=-\delta_{i,j}$. D'o\`u l'on d\'eduit
l'assertion voulue sur la signature.
\ste (iii)
On utilisera le fait suivant :\ste
{\em Si $\rho$ est une repr\'esentation complexe d'une
alg\`ebre de Lie simple complexe $\s$ dans un espace de dimension
finie $V$, on a :$$tr(\rho(X)\rho(Y))=qK_{\s}(X,Y)$$
o\`u $q$ est un nombre rationel positif. De plus $q$ est nul si et
seulement si $\rho$ est triviale.}\ste
L'existence d'un coefficient de proportionnalit\'e $q$ est claire ,
car
la forme de Killing est la seule forme $\C$-bilin\'eaire sur $\s$, \`a
un scalaire multiplicatif pr\`es. On se ram\`ene, pour l'\'etude de $q$,
au cas o\`u  $\rho$ est simple. On consid\`ere, sur $V$, un produit
scalaire invariant par une forme r\'eelle compacte, $\k$, de $\s$. Si
$X\in
\k$, $\rho(X)$ est antihermitien et :
$$tr(\rho(X) \rho(X))=-tr(\rho(X) \rho(X)^{*})\leq 0$$
cette trace \'etant nulle seulement si $\rho(X)$ est nul. On en
d\'eduit que $q>0$ si $\rho$ est non triviale. Puis on prend un
\'el\'ement
non nul d'une sous-alg\`ebre de Cartan $\j$ de $\s$, sur lequel tous
les poids entiers de $\j$   sont entiers. On en d\'eduit que
$K_{\g}(X,X)$ et
$tr(\rho(X)
\rho(X))$ sont des entiers, le premier nombre \'etant non nul car
\'egal \`a la somme sur toutes les racines, $\alpha$, de
$(\alpha(X))^{2}$. La rationalit\'e de $q$ en r\'esulte.\qed
\begin{defi}\ste Si $\g$ est une alg\`ebre de Lie r\'eductive complexe,
une forme
$\R$-bilin\'eaire\ste sym\'etrique sur $\g$ et invariante par  $\g$ est
dite forme de Manin si et seulement si elle est de signature
$(dim_{\C}\g,dim_{\C}\g)$. Une forme de Manin est dite forme
sp\'eciale si sa restriction \`a toute sous-alg\`ebre
de Lie complexe semi-simple est non d\'eg\'en\'er\'ee.
\end{defi}

\begin{lem}
\ste (i) Une
forme
$\R$-bilin\'eaire sym\'etrique $\g$-invariante sur $\g$
est sp\'eciale
si et seulement si sa restriction \`a $\g^{der}$ est sp\'eciale et si
sa restriction au centre, $\z$, de $\g$ est de signature
$(dim_{\C}\z,dim_{\C}\z)$.\ste
(ii) La restriction d'une forme sp\'eciale \`a une sous-alg\`ebre de Lie
semisimple complexe de $\g$ est sp\'eciale.\ste
 (iii) La restriction
d'une forme sp\'eciale au centralisateur d'un \'el\'ement semi-simple de
$\g$,  dont l'image par la repr\'esentation adjointe de $\g$ n'a que
des  valeurs propres  r\'eelles,
 est sp\'eciale.\ste (iv) Si $\g$ est semi-simple, et
$B=B_{\lambda}^{\g}$, o\`u
$\lambda=(\lambda_{1},\dots,
\lambda_{n})\in\C^{n}$ v\'erifie :
\begin{equation}Si\>\> \sum_{i=1,\dots, n}q_{i}\lambda_{i}=0\>\>
avec\>\> q_{i}\in
\Q^{+},\>\> alors \>\>les \>\>q_{i}\>\> sont\>\> tous\>\>
nuls\end{equation}
alors $B$ est sp\'eciale.
On note que (1.1) est satisfait d\`es que les $\lambda_{i}$ sont
ind\'ependants sur $\Q$, ou bien tous strictement positifs. \ste
(v) Si $\g$ est simple, toute forme
$\R$-bilin\'eaire sym\'etrique $\g$-invariante sur $\g$ est sp\'eciale.
\end{lem}
{\em D\'emonstration :} (i) r\'esulte du fait
que toute sous-alg\`ebre semi-simple de $\g$ est
contenue dans $\g^{der}$ et que, pour toute forme
$\R$-bilin\'eaire sym\'etrique $\g$-invariante sur
$\g$, le centre de $\g$ et $\g^{der}$ sont
orthogonaux.\ste
(ii) est clair. \ste
Montrons  (iii). Comme le centralisateur d'un \'el\'ement de $\g$
est la somme de son intersection avec $\g^{der}$
et de celle avec le centre, on se r\'eduit ais\'ement,
gr\^{a}ce \`a (i) au cas o\`u $\g$ est semi-simple, ce que
l'on suppose dans la suite. Soit $X$ un \'el\'ement semi-simple de
$\g$ tel que $ad\>X$ n'a que des valeurs propres r\'eelles,
soit $\l$ son centralisateur et $\c$ le centre de $\l$. D'apr\`es
(i) et (ii), il suffit de voir que la restriction d'une forme
sp\'eciale \`a $\c$ est de signature $(dim_{\C}\c,dim_{\C}\c)$. Soit
$\j$ une sous-alg\`ebre de Cartan   de $\g$ contenant $X$. Alors
$\c$ est \'egal \`a l'intersection des noyaux des racines de $\j$
dans $\g$ s'annulant sur $X$. Cela montre que $\c$ est la
somme de ses intersections avec les id\'eaux simples de $\g$.
Il suffit alors de prouver notre assertion sur la signature
dans le cas o\`u $\g$ est simple. Alors $B=B_{\lambda}$, avec
$\lambda$ non nul. Soit $\j_{\R}$ l'espace form\'e des \'el\'ements de
$\j$ sur lesquels toutes les racines de $\j$ dans $\g$ sont
r\'eelles, qui est une forme r\'eelle de $\j$. Il est clair que
$\c$ est la somme directe de $\c_{\R}:=\c\cap\j_{\R}$ avec
$i\c_{\R}$. On fixe
 une
base orthonorm\'ee,  $X_{1},\dots, X_{l}$,  de $\c_{\R}$, pour la
forme de Killing de $\g$. Celle-ci existe car la forme de Killing
est d\'efinie  positive sur
$\j_{\R}$. On choisit une racine carr\'ee $\mu$ de $i\lambda^{-1}$.
Onpose
$Y_{i}=\mu X_{i}$, $Z_{i}=i\mu X_{i}$. Alors
 $B_{\lambda}(Y_{i},Z_{j})=0$,
$B_{\lambda}(Y_{i},Y_{j})=\delta_{i,j}$,
$B_{\lambda}(Z_{i},Z_{j})=-\delta_{i,j}$. D'o\`u l'on d\'eduit
l'assertion voulue sur la signature, ce qui prouve (iii).\ste
(iv) est une cons\'equence imm\'ediate du Lemme 1 et (v) est un cas
particulier de (iv)\qed

\begin{cor}
(i) La restriction,  d'une forme
$\R$-bilin\'eaire, sym\'etrique, $\g$-invariante sur $\g$, et non
d\'eg\'en\'er\'ee,
$B$ , au
 centralisateur, $\l$ d'un \'el\'ement semi-simple de
$\g$,  dont l'image par la repr\'esentation adjointe de $\g$ n'a que
des  valeurs propres  r\'eelles, est non d\'eg\'en\'er\'ee. Il en va de
m\^{e}me de sa restriction \`a $\l^{der}$ et au centre $\a$ de $\l$
.\ste (ii) Si
$B$ est une forme de Manin sur $\g$, sa restriction \`a
$\l$ (resp. $\l^{der}$, $\a$) est une forme de Manin sur $\l$
(resp. $\l^{der}$, $\a$).\ste
(iii) Une forme bilin\'eaire sym\'etrique, $\g$-invariante est une
forme de Manin si et seulement si sa restriction \`a $\g^{der}$ et sa
restriction au centre de $\g$ sont des formes de Manin.
\end{cor}
\dem
L'alg\`ebre de Lie $\l$ est la somme du centre $\z$ de $\g$ avec la
somme de ses intersections avec les id\'eaux simples de $\g$. Ces
sous-alg\`ebres sont deux \`a deux orthogonales pour $B$, d'apr\`es le
Lemme 1. La restriction de $B$ \`a $\z$ est non d\'eg\'en\'er\'ee, car
$\g^{der}$ et $\z$ sont orthogonaux. De plus la restriction de
$B$ \`a l'intersection de $\l$ avec un id\'eal simple de $\g$ est non
d\'eg\'en\'er\'ee, d'apr\`es le Lemme 2 (iii), appliqu\'e \`a cet id\'eal
simple. D'o\`u l'on d\'eduit que la restriction de $B$ \`a $\l$ est non
d\'eg\'en\'er\'ee, ce qui implique le m\^{e}me fait pour sa restriction \`a
$\l^{der}$ et $\a$.
\ste Si $B$ est de signature $(dim_{\C}\g, dim_{\C}\g) $, sa
restriction \`a $\z$ est de signature $(dim_{\C}\z, dim_{\C}\z) $.
Alors, les assertions sur la signature se d\'emontre comme ci-dessus,
gr\^{a}ce au Lemme 2 (ii), (v), apliqu\'e aux id\'eaux simples de
$\g$. D'o\`u (ii). \ste
La partie si de (iii) est claire. La partie seulement si r\'esulte
de (ii) \qed \ste
On rappelle que le radical  d'une alg\`ebre de Lie,
$\g$ est son plus grand id\'eal r\'esoluble,
et que son radical nilpotent,
  est le plus grand id\'eal, dont les \'el\'ements sont repr\'esent\'es
par des endomorphismes nilpotents dans chaque repr\'esentation de
dimension finie de $\g$.
 Suivant Bourbaki, on
appelle sous-alg\`ebre de Levi d'une alg\`ebre de Lie, toute
sous-alg\`ebre semi-simple suppl\'e-\ste mentaire du radical. Si $\g$
est une alg\`ebre de Lie semi-simple complexe on appelle d\'ecomposition
de Langlands d'une sous-alg\`ebre parabolique $\p$ de $\g$ une
d\'ecomposition de la forme $\p=\l\oplus\n$, o\`u $\n$ est le
radical nilpotent et $\l$ est une sous-alg\`ebre de Lie complexe de
$\g$, r\'eductive dans
$\g$. \ste
Rassemblons dans un Lemme quelques propri\'et\'es des
d\'ecompositions de Langlands d'une sous-alg\`ebre parabolique.
\ste

\begin{lem}
Soit une $\p$ sous-alg\`ebre parabolique  de $\g$, $\n$ son
radical nilpotent.\ste (i) Si $\j$ est une sous-alg\`ebre de Cartan
de $\p$, c'est une sous-alg\`ebre de Cartan
de $\g$  , et il existe une seule d\'ecomposition de
Langlands  de $\p$,
$\p=\l\oplus\n$, telle que $\j$ soit contenue dans $\l$.
\ste (ii) Si  $\j$, $\j'$  sont  deux sous-alg\`ebres de Cartan de
$\g$, contenues dans $\p$, elles sont conjugu\'ees par un \'el\'ement de
$P$, qui conjugue les alg\`ebres $\l$ et $\l'$
correspondantes.\ste  (iii) Si $\p=\l\oplus\n$ est une
d\'ecomposition de Langlands de
$\p$, toute sous-alg\`ebre de Cartan
de $\p$ est une  sous-alg\`ebre de Cartan
de $\g$.
\end{lem}
\dem Revenant \`a la d\'efinition des sous-alg\`ebres paraboliques (cf
[Bou], Ch. VIII, Paragraphe 3.4, D\'efinition 2, par exemple), on
voit qu'il existe   une sous-alg\`ebre de Cartan
de $\g$, $\j_{1}$, et une d\'ecomposition de Langlands de $\p$,
$\p=\l_{1}\oplus \n$, avec $\j_{1}$ contenue dans $\l_{1}$, telle
que $\l_{1}$ soit la somme des sous-espaces poids de  $\j_{1}$
dans $\p$, qui ne rencontre pas $\n$. En particulier les poids de
$\j_{1}$ dans $\l_{1}\approx \p/\n$ sont distincts de ceux dans
$\n$. Si
$\p=\l'_{1}\oplus \n$ est une d\'ecomposition de Langlands de $\p$,
avec $\j_{1} $ contenue dans $\l'_{1}$, $\l'_{1}$ est somme des
sous-espaces poids de $\j_{1}$ dans $\l'_{1}\approx \g/\n$. D'o\`u
l'\'egalit\'e de  $\l_{1}$ et $\l'_{1}$. Ceci assure l'unicit\'e de
$\l$ pour $\j=\j_{1}$. Maintenant toutes les sous-alg\`ebres de
Cartan de
$\p$ sont conjugu\'ees \`a $\j_{1}$, par un \'el\'ement de $P$ (cf.
 [Bour],
Paragraphe 3.3, Th\'eor\`eme 1). On en d\'eduit (i) par transport de
structure, et
(ii) r\'esulte de la preuve de (i).\ste
Montrons (iii). Si $\p=\l\oplus\n$ est une d\'ecomposition de
Langlands de
$\p$, l'action du centre de $\l$ sur $\g$ est semi-simple, puisque
$\l$ est r\'eductive dans $\g$. Cela implique que, si $\j$ une une
sous-alg\`ebre de Cartan de $\l$, $\j$ est ab\'elienne et  form\'ee
d'\'el\'ements semi-simples de $\g$. Mais $\l$ est isomorphe \`a $\p/\n$,
qui est une alg\`ebre r\'eductive de m\^{e}me rang que $\g$. Pour des
raisons de dimension, on voit que $\j$ est une sous-alg\`ebre de
Cartan de $\g$.\qed
 \begin{defi} On appelle af-involution  (a pour antilin\'eaire, f
pour flip) d'une alg\`ebre de Lie semi-simple complexe $\m$, tout
automorphisme involutif, $\R$-lin\'eaire, $\sigma$, de $\m$, pour
lequel  il existe :\ste
1) un id\'eal ${\tilde  \m_{0}}$ de  $\m$, et une forme r\'eelle,
${\tilde \h}_{0}$ of $\m_{0}.$\ste   2) des id\'eaux simples de
$\m$, $\m'_{j}$,
$\m''_{j}$, $j=1,\dots,p$.\ste
3) un  isomorphisme d'alg\`ebres de Lie complexes, $\tau_{j}$,
entre
$\m'_{j}$ et
$\m''_{j}$, $j=1,\dots,p$, \ste
tel que, notant  $ \h_{j}:=\{(X, \tau_{j}(X)\vert X
\in
\m'_{j}\}$, et notant  $\h$, l'ensemble des points fixes de
$\sigma$, on ait :

$$\m={\tilde \m_{0}}\oplus (\oplus_{ j=1,\dots,p} (\m'_{j}\oplus
\m''_{j}))$$
$$\h= {\tilde \h}_{0}\oplus (\oplus_{ j=1,\dots,p} \h_{j}) $$
\end{defi}
Il est bon de remarquer qu'un automorphisme involutif de $\m$ est
caract\'eris\'e par son espace de points fixes, car l'espace des
\'el\'ements anti-invariants est juste l'orthogonal de celui-ci, pour
la forme de Killing de $\m$ regard\'ee comme r\'eelle.\ste
D\'ebutons par quelques propri\'et\'es \'el\'ementaires.
\begin{lem}
Soit $\h$ une forme r\'eelle simple d'une alg\`ebre de Lie semi-simple
complexe $\s$. \ste (i) L'alg\`ebre $\s$ n'est pas simple, si et
seulement si
$\h$ admet une structure complexe \ste  (ii) Dans ce cas, $\s$ est
le produit   de deux id\'eaux simples, $\s_{1}$, $\s_{2}$, isomorphes
\`a
$\h$. \ste (iii) Toujours dans ce cas, il existe un isomorphisme
antilin\'eaire,
$\tau$, entre les alg\`ebres de Lie
$\s_{1}$, $\s_{2}$, regard\'ees comme r\'eelles, tel que
: $$\h=\{(X,\tau (X))\vert X\in \s_{1}\}$$
\end{lem}
\dem Les  points (i) et (ii) sont  bien connus. Montons (iii). Comme
$\h$ est une forme r\'eelle de $\s_{1}\oplus \s_{2}$, la projection de
$\h$ sur chacun des deux facteurs est non nulle, donc induit un
isomorphisme $\R$-lin\'eaire  de $\h$ avec chacun de ces facteurs.
Il en r\'esulte que $\h$ a la forme indiqu\'ee, mais on sait seulement
que $\tau $ est $\R$-lin\'eaire. Mais alors $\h$ appara\^{i}t comme
l'ensemble des points fixes de l'automorphisme involutif
$\R$-lin\'eaire de
$\s$ d\'efini par : $$(X,Y)\mapsto (\tau ^{-1}(Y),\tau (X)), X\in
\s_{1}, Y\in \s_{2}
$$
D'apr\`es la remarque qui pr\'ec\`ede le Lemme, cette involution doit
\^{e}tre \'egale \`a la conjugaison par rapport \`a  $\h$, donc elle est
antilin\'eaire. Ceci implique l'antilin\'earit\'e de $\tau$. \qed

\begin{lem}
On se donne une af-involution, $\sigma$, d'une alg\`ebre de Lie
semi-simple de $\m$. Les id\'eaux simples de $\m$ sont permut\'es par
$\sigma$. On note
$\m_{j}$,
$j= 1,
\dots,r$, les id\'eaux simples de $\m$. On d\'efinit une involution
$\theta$ de
$\{1,
\dots,r\}$ caract\'eris\'ee par  : $\sigma (\m_{j})=\m_{\theta (j)},
\>\>j =1,
\dots,r$. Pour $j =1,
\dots,r$, l'une des propri\'et\'es suivantes est v\'erifi\'ee :\ste
1)  $\theta (j)=j$ et la restriction de $\sigma $ \`a $\m_{j}$ est
un automorphisme antilin\'eaire de $\m_{j}$. \ste
2) $\theta (j)\not =j$ et la restriction de $\sigma $ \`a $\m_{j}$ est
un isomorphisme antilin\'eaire de $\m_{j}$ sur $\m_{\theta
(j)}$. \ste   3) $\theta (j)\not =j$ et la restriction de $\sigma $
\`a
$\m_{j}$ est un isomorphisme $\C$-lin\'eaire de $\m_{j}$ sur
$\m_{\theta (j)}$.\ste
On est dans le cas 1) ou 2) si et seulement si $\m_{j}$ est
contenu dans l'id\'eal ${\tilde \m}_{0}$ de la d\'efinition des
af-involutions.
\end{lem}
\dem
En effet, soit $\h_{p+l}$, $l=1, \dots , q$, les id\'eaux simples de
${\tilde \h_{0}}$. Comme ${\tilde \h_{0}}$ est une forme r\'eelle de
${\tilde \m_{0}}$, ${\tilde \m_{0}}$ est la somme directe des
$\h_{l}+i\h_{l}$, qui sont en outre des id\'eaux. Si
$\h_{l}+i\h_{l}$ est simple, c'est un id\'eal simple de $\m$ et on
est dans le cas 1). Sinon $\h_{l}+i\h_{l}$ est le produit de deux
id\'eaux simples $\m_{j}\times \m'_{j}$ et l'on est dans le cas
2). \ste On traite de m\^{e}me le cas o\`u
$\m_{l}$ est
\'egal
\`a l'un des
$\m'_{j}$,
$\m''_{j}$,
$j=1, \dots \>\>, p$, en remarquant que $\h_{j}$ \ste est l'ensemble
des points fixes de l'involution $\C$-lin\'eaire de $(\m'_{j},
\m''_{j})$  donn\'ee par :
\beq(X,Y)\mapsto (\tau ^{-1}(Y),\tau (X)), X\in
\m'_{j}, Y\in \m'' _{j}
\eeq
\qed
\begin{lem}
Tout isomorphisme $\R$-lin\'eaire entre deux alg\`ebres de Lie
simples complexes est soit $\C$-lin\'eaire, soit antilin\'eaire.
\end{lem}
\dem
Deux alg\`ebres de Lie semi-simples complexes qui sont isomorphes
comme alg\`ebres r\'eelles, sont isomorphes comme alg\`ebre de Lie
complexes. En effet, leurs syst\`emes de racines restreintes sont
alors isomorphes. Mais chacun de  ceux-ci est  isomorphe au syst\`eme
de racine  associ\'e \`a une \sac. D'o\`u l'assertion. Ceci implique
que l'on peut se limiter, pour prouver le Lemme, aux automorphismes
$\R$-lin\'eaire  d'une alg\`ebre de Lie simple complexe, $\g$.
Consid\'erant la conjugaison par rapport \`a une forme r\'eelle de
$\g$, $X\mapsto {\overline X}$, l'alg\`ebre de Lie
$\g':=\{(X,{\overline X)}\vert X\in \g\}$ est une forme r\'eelle de
$\g\times \g$ isomorphe \`a $\g$. Alors tout automorphisme, $\sigma$,
$\R$-lin\'eaire de
$\g$, d\'efinit,  par transport de  structure, un automorphisme de
$\g'$, $\sigma'$, qui poss\`ede un unique prolongement $\C$-lin\'eaire
\`a
$\g\times \g$, $\sigma''$ . Il existe deux automorphismes
$\C$-lin\'eaires de
$\g$, $\tau $ et $\sigma$, tels que $\sigma''$ v\'erifie :
$$\sigma''(X,Y)=(\tau(X), \theta(Y)) \>\>ou \>\>bien\>\> (\tau(Y),
\theta(X)), (X,Y)\in \g\times \g$$
Ecrivant la d\'efinition de $\sigma'$, la stabilit\'e de
$\g'$, par
$\sigma'' $ implique que $\sigma$ est $\C$-lin\'eaire dans le premier
cas et antilin\'eaire dans le second.\qed \ste
{\bf Corollaire }
{\em Une involution d'une alg\`ebre de Lie semi-simple complexe est
une af-involution si et seulement si sa restriction \`a tout id\'eal
simple qu'elle laisse invariant est antilin\'eaire.}\ste
\dem En effet, d'apr\`es le Lemme 5, il suffit de voir que tout
automorphisme du produit de deux alg\`ebres de Lie simples complexes
$\s_{1}$, $\s_{2}$, permutant les facteurs, est soit $\C$-lin\'eaire,
soit antilin\'eaire. Mais un tel automorphisme est de la forme :
$$(X,Y)\mapsto (\tau ^{-1}(Y),\tau (X)), X\in
\s_{1}, Y\in \s_{2},
$$
o\`u $\tau$ est un isomorphisme $\R$-lin\'eaire entre $\s_{1}$, $\s_{2}$
. On conclut gr\^{a}ce au Lemme pr\'ec\'edent.\qed
Comme on l'a indiqu\'e dans l'introduction, le Th\'eor\`eme suivant
g\'en\'eralise des r\'esultats d'E. Karolinsky (cf. [K1],
Th\'eor\`eme 3 (i) et [K3] Proposition 3.1)
\begin{theo} Soit $\g$ une alg\`ebre de Lie r\'eductive complexe et
$B$ une forme de Manin. Soit $\i$ une sous-alg\`ebre de Lie r\'eelle de
$\g$, Lagrangienne pour
$B$.\ste On a les propri\'et\'es suivantes  :\ste (i) Si l'on note $\p$
le normalisateur dans $\g$ du radical nilpotent, $\n$, de $\i$, $\p$
est une sous-alg\`ebre parabolique de $\g$, contenant $\i$, de radical
nilpotent
$\n$.\ste (ii) Soit $\p=\l\oplus \n$ une d\'ecomposition de
Langlands de $\p$, $\a$ le centre de $\l$ et  $\m$ son id\'eal
d\'eriv\'e. On note   $\h$ l'intersection de $\i$ et $\m$. Elle est
isotrope  pour
$B$.
De plus $\h$ est l'espace des points fixes d'une af-involution de
$\m$. \ste Si $B$ est sp\'eciale, celle-ci est antilin\'eaire et $\h$
est une forme r\'eelle de $\m$.
\ste(iii) L'intersection $\i_{\a}$ de $\a$ et $\i$ est
Lagrangienne  pour la restriction de $B$ \`a $\a$.
\ste (iv) On a $\i=\h\oplus\i_{\a}\oplus \n$.\ste\ste
 R\'eciproquement si une sous-alg\`ebre de Lie r\'eelle, $\i$, de
$\g$ est de la forme ci-dessus, elle est Lagrangienne pour $B$.
 On dit alors que {\bf  $\i$ est sous $\p$}.
\end{theo}

{\em D\'ebut de la d\'emonstration du Th\'eor\`eme 1 :}
 Soit $\i$  une sous-alg\`ebre de
Lie r\'eelle de $\g$, isotrope pour $B$, dont la dimension est
\'egale \`a la dimension complexe de $\g$. On note $\r$ son
radical et on pose :  \beq
\n:=\{ X\in \r\cap\g^{der}\vert ad_{\g}(X)\>\> est\>\>
nilpotent\}\eeq
Soit  $\h$ une sous-alg\`ebre
de Levi de $\i$.
\begin{lem}
L'ensemble $\n$ est un id\'eal de $\i$ et $[\i,\r]$ est contenu
dans $\n$.
\end{lem}
{\em D\'emonstration :} Montrons que $\n$ est un id\'eal de $\r$
contenant $[\r, \r]$. En effet, comme $\r$ est r\'esoluble, dans
une base, sur $\C$,  bien choisie de $\g$, les $ad_{\g}(X)$,
$X\in \r$ s'\'ecrivent sous forme de matrices triangulaires
sup\'erieures. Pour $X\in \r$, les entr\'ees de la diagonale de
cette matrice sont not\'ees $\lambda_{1}(X),\dots,\lambda_{p}(X)
$, o\`u les $\lambda _{i}$ sont des caract\`eres de $\r$. Alors
$\n$ est l'intersection des noyaux de ces caract\`eres avec
$\g^{der}$. Donc  $\n$ est  un id\'eal de $\r$ contenant $[\r,
\r]$.\ste Si $\f$ est une sous-alg\`ebre de Cartan de $\h$,
$\f\oplus\r$ est encore une alg\`ebre de Lie r\'esoluble car
$[\f,\f]=\{0\}$ et $[\f, \r]\subset \r$. Un argument similaire \`a
celui ci-dessus montre que $[\f, \r]$ est contenu dans $\n$. La
r\'eunion de toutes les sous-alg\`ebres de Cartan de $\h$ est
dense dans $\h$, d'apr\`es la densit\'e des \'el\'ements r\'eguliers
(cf.[Bou], Ch. VII, Paragraphe 2.2 et Paragraphe 2.3, Th\'eor\`eme 1).
 Par continuit\'e
et densit\'e, on en d\'eduit que $[\h,\r]\subset \n$.\qed

\begin{lem}Soit $k$ un entier compris entre 0 et la dimension
r\'eelle de $\r / \n$. Il existe un sous-espace r\'eel, ab\'elien,
$\a_{k}$,  de $\r$, de dimension $k$, tel que :\ste
(i) $\a_{k}\cap\n=\{ 0\}.$
\ste(ii) $\a_{k} $ est form\'e d'\'el\'ements semi-simples de $\g$.
\ste (iii) $\a_{k} $ et $\h$ commutent. \end{lem}  {\em
D\'emonstration : } On proc\`ede par r\'ecurrence sur $k$. Si $k=0$,
le Lemme est clair. Supposons le d\'emontr\'e pour
$k<dim_{\R}(\r/\n)$ et montrons le au rang $k+1$. Alors $\h\oplus
\a_{k}$ est une alg\`ebre de Lie r\'eductive dans $\g$, regard\'ee
comme r\'eelle. D'autre part, comme $\n$ contient $[\i, \r]$
d'apr\`es
le Lemme pr\'ec\'edent, on voit que
$\i$  et donc  $\h\oplus \a_{k}$ agit trivialement sur $\r/\n$.
Ceci
implique que : $$[\h\oplus \a_{k},\a_{k}\oplus\n]\subset \n$$
Donc, $ \a_{k}\oplus\n$ est un $(\h\oplus \a_{k})$-sous-module
de $\r$, qui admet un suppl\'ementaire dans $\r$ commutant \`a
$\h\oplus \a_{k}$, puisque
$\h\oplus \a_{k}$ est r\'eductive dans $\g$ et que le quotient
$\r/\a_{k}\oplus\n$ est un $\h\oplus \a_{k}$-module trivial. \ste
On choisit un \'el\'ement non nul de ce suppl\'ementaire, $X$. Alors :
\begin{equation}X\in \r, \>\>X\notin \h\oplus \a_{k}, \>\> et
\>\>[X,\h\oplus \a_{k}]=\{0\} \end{equation} On \'ecrit
$X=X_{s}+X_{n}$, o\`u $X_{n}$ est un \'el\'ement de $\g^{der}$,
$X_{s}$ est un \'el\'ement de $\g$  commutant \`a $X_{n}$
 tels que
$ad_{\g}X_{s}$ est semi-simple et $ad_{\g}X_{n}$est nilpotent.
On sait qu'alors $ad_{\g}X_{s}$, $ad_{\g}X_{n }$ sont des
polyn\^omes en $ad_{\g}X$. Joint \`a (1.2), cela implique :
\begin{equation}[X_{s},\h\oplus
\a_{k}]=\{0\},\>\>[X_{n},\h\oplus
\a_{k}]=\{0\}
\end{equation}
Montrons que $X_{n}$ appartient \`a $\i$. Soit $\j$ une sous-alg\`ebre
de Cartan de $\h$. Alors $\j\oplus \r$ est r\'esoluble. On peut donc
choisir une base de $\g$ dans laquelle les
$ad_{\g}Y$, $Y\in \j\oplus \r$, sont repr\'esent\'es par des
matrices triangulaires sup\'erieures. On peut choisir cette base
de sorte qu'elle soit la r\'eunion de bases des id\'eaux simples de
$\g$ avec une base du centre de $\g$, ce que l'on fait dans la
suite. Comme  $ad_{\g}X_{n }$ est un polyn\^ome en $ad_{\g}X$, et
que $X\in \r$, il est repr\'esent\'e dans cette base par une
matrice triangulaire sup\'erieure. Comme cet endomorphisme est
nilpotent, sa diagonale est nulle. On en d\'eduit que, pour tout
$Y\in \j\oplus\r$,  les composantes de $X_{n}$ et $Y$ dans les
id\'eaux simples  de $\g$ sont deux \`a deux
orthogonales pour la forme de Killing de $\g$.
Alors, il r\'esulte de l'orthogonalit\'e, pour $B$,  du centre de
$\g$  \`a $\g^{der}$ et du Lemme 1 (i), que :
$$B(X_{n},Y)=0, \>\> Y\in \j\oplus \r$$
En utilisant la densit\'e dans $\h$ de la r\'eunion de ses
sous-alg\`ebres de Cartan, on en d\'eduit que $X_{n}$ est
orthogonal \`a $\i$ pour $B$. Mais $\i$ est un  sous-espace
isotrope pour $B$ de dimension maximale, d'apr\`es la d\'efinition
des formes de Manin. Donc $X_{n}$ est \'el\'ement de $\i$ comme
d\'esir\'e.\ste
Ecrivons $X_{n}=H+R$ avec  $H\in \h$, $R\in  \r$. Comme
$X_{n}$ commute \`a $\h$, d'apr\`es (1.3) et que $[\h, \r]\subset
\r$, on voit que $H$ commute \`a $\h$. Donc  $H$ est nul puisque
$\h$ est semi-simple. Finalement $X_{n}\in \r$, et en fait
$X_{n}\in\n$, d'apr\`es la d\'efinition de $\n$. Comme $X$
appartient \`a un suppl\'ementaire de $\a_{k}+\n$ dans $\r$ et que
$X=X_{s}+X_{n}$, on a :
$$X_{s}\in \r, X_{s}\notin \a_{k}+\n$$
On pose $\a_{k+1}=\a_{k}+\R X_{s}$. D'apr\`es (1.3) et la
semi-simplicit\'e de $ad_{\g}X_{s}$, $\a_{k+1}$ v\'erifie les
propri\'et\'es voulues.\qed
\ste{\em Suite de la d\'emonstration du Th\'eor\`eme 1:}\ste On pose
$\i_{\a}:= \a_{p}$, avec $p=dim_{\R}\r/\n$, de sorte que
$\i=\h\oplus \i_{\a}\oplus \n$, o\`u $\i_{\a}$ est form\'e d'\'el\'ements
semi-simples de $\g$ avec : $$[\h,
\i_{\a}]=\{0\},\>\>\r=\i_{\a}\oplus \n$$ Comme $\i_{\a}\oplus \n$ est
r\'esoluble, il existe une sous-alg\`ebre de Borel, $\b$, de $\g$,
contenant $\i_{\a}\oplus \n$. A noter que $\n$ est contenue dans le
radical nilpotent, $\v$, de $\b$, d'apr\`es la d\'efinition de $\n$
et les propri\'et\'es du radical nilpotent d'une sous-alg\`ebre de
Borel. Montrons que
$\i_{\a}$ est contenue dans une sous-alg\`ebre de Cartan de $\g$,
contenue dans
$\b$. En effet, d'apr\`es  [Bor], Proposition 11.15, la sous-alg\`ebre
de Borel $\b$ contenant $\i_{\a}$, elle contient une sous-alg\`ebre
de Borel du centralisateur $\l$ de $\i_{\a}$ dans
$\g$. Celle-ci contient une sous-alg\`ebre de Cartan $\j$ de
$\l$. Celle-ci est aussi une sous-alg\`ebre de Cartan
de $\g$ contenant $\i_{\a}$ (cf. [Bou], Ch. VII, Paragraphe 2.3,
Proposition 10). \ste Soit $\u$ la somme  des sous
espaces poids de $\i_{\a}$, dans $\v$, , pour des poids non nuls.
Alors
$\p:=\l\oplus\u$ est une sous alg\`ebre parabolique de $\g$,
contenant $\b$. Comme $\l$ est r\'eductive,
$\m:=\l^{der}$ est semi-simple et le radical de $\p$ est \'egal \`a la
somme du centre $\a$ de $\l$ avec $\u$. La d\'efinition de $\u$
montre que $[\p,\p]=\m\oplus \u$, donc le radical nilpotent de $\p$
est \'egal \`a $\u$ (cf. [Bou], Ch. I, Paragraphe 5.3, Th\'eor\`eme 1).
 Comme
 $\i=\h\oplus \i_{\a}\oplus \n$, que $\i_{\a}\oplus \n$ est
contenu dans $\b$ et  que $\h$ est contenu dans $\l$, on a :
$$\i\subset\p$$
Or $\p$ (resp. $\u$) est la somme de ses intersections
$\p_{i}$ (resp. $\u_{i}$) avec les id\'eaux simples $\g_{i}$ de
$\g$. Comme $\p_{i}$ est orthogonal \`a $\u_{i}$ pour la forme de
Killing de $\g_{i}$, on en d\'eduit que $\u$ est orthogonal \`a $\p$
pour $B$ (cf. Lemme 1 (i)). Comme $\i$ est un sous-espace
isotrope pour $B$, de dimension maximale et contenu dans $\p$,
$\u$ est inclus dans $\i$. Il r\'esulte alors de la d\'efinition de
$\n$, que $\u$ est contenu dans $\n$. Par suite, on a  :
\begin{equation}
\i=\u\oplus
(\i\cap\l),\>\>\n=\u\oplus(\n\cap\l)
\end{equation}
 Remarquons que $\a$  contient
$\i_{\a}$.
 On a $\v=\u\oplus
(\v\cap\m)$. Comme $\n\subset \v$ et $\u\subset\n$, on en d\'eduit
que $\n=\u\oplus (\n \cap\m)$. On d\'eduit alors de (1.6) que :
$\n\cap \l=\n\cap\m$. Finalement, on a :
$$\i=\h\oplus \i_{\a}\oplus (\n\cap\m)\oplus \u$$
Alors $\i':=\i\cap\m $ est \'egal \`a $\h\oplus(\n\cap\m)$. C' est une
sous-alg\`ebre isotrope de $\m$ pour la restriction de $B$ \`a $\m$,
donc, d'apr\`es le Corollaire du Lemme 2 (ii), de dimension r\'eelle
inf\'erieure ou \'egale \`a la dimension complexe de $\m$. \ste
De m\^{e}me, $\i_{\a}$ est un sous espace isotrope de $\a$ pour la
restriction de $B$ \`a $\a$. D'apr\`es le Corollaire du Lemme 2, la
restriction de $B$ \`a $\a$ est de signature $(dim_{\C}\a,
dim_{\C}\a)$. Il en r\'esulte que la
dimension de $\i_{\a}$ est inf\'erieure ou \'egale \`a $dim_{\C}\a$.
Mais $dim_{\C}\g=dim_{\C}\m+dim_{\C}\a+dim_{\R}\u$. Comme
$dim_{\R}\i=dim_{\C}\g$, on d\'eduit de ce qui pr\'ec\`ede que l'on a
: \begin{equation}dim_{\R} \i'=dim_{\C}\m,\>\>
dim_{\R}\i_{\a}=dim_{\C}\a
\end{equation}
\begin{lem}L'alg\`ebre de Lie  $\n':=\n\cap\m$ est r\'eduite  \`a z\'ero
et $\h$ \`a la forme indiqu\'ee dans le Th\'eor\`eme.
\end{lem}{\em D\'emonstration :}
Si
$\f$ est une sous-alg\`ebre de Cartan de
$\h$,
$\f+\n'+\i\n'$ est une sous-alg\`ebre de Lie r\'eelle et r\'esoluble de
$\m$. On peut donc choisir une base de $\m$, r\'eunion de bases
des id\'eaux simples de $\m$ , telle que, pour
tout $X\in\f+\n'+\i\n'$, $ad_{\m}X$ soit repr\'esent\'e, dans cette
base, par une matrice triangulaire sup\'erieure. De plus, si $X$
est \'el\'ement de $\n'+\i\n'$, les \'el\'ements diagonaux de cette
matrice sont nulles. On voit, gr\^{a}ce au Lemme 1 (i), que
$\n'+i\n'$ est orthogonal \`a $\f+\n'$, pour la restriction, $B'$, de
$B$ \`a $\m$. Ceci \'etant vrai pour tout $\f$, $\n'+i\n'$
est  orthogonal \`a
$\i'$ ($=\h+\n'$), pour $B'$. Mais $B'$ est non d\'eg\'en\'er\'ee, d'apr\`es
le Corollaire du Lemme 2, donc, d'apr\`es le Lemme 1 (ii) et (1.7),
$\i'$ est un sous-espace isotrope, pour
$B'$, de dimension maximale de $\m$. Il en r\'esulte que
$\n'+i\n'$ est contenu dans $\i'$. Mais $\n'+i\n'$ est aussi
  contenu dans $ \v\subset \g^{der}$. Finalement $\n'+i\n'$ est
 contenu dans l'intersection  de $\i'$ avec  $\n$, d'apr\`es
la d\'efinition de celui-ci. Mais, comme $\i'$ est contenu dans $\m$,
$\i'\cap\n =\n'$. Alors on a :
$\n'+i\n'\subset
\n'$, c'est \`a dire que $\n'$ est un sous-espace vectoriel
complexe de $\g$. \ste Soit $\h_{j}$, $j=1,\dots,r$, les id\'eaux
simples de $\h$. Comme  $\h_{j}\cap i \h_{j}$ est un id\'eal de
l'alg\`ebre de Lie simple r\'eelle $\h_{j}$, il y a deux possiblit\'es pour
$\h_{j}$. Ou bien $\h_{j}\cap i \h_{j}=\{0\}$, et alors $\h_{j}+
i \h_{j}$ est une alg\`ebre de Lie semi-simple complexe dont
$\h_{j}$ est une forme r\'eelle. Ou bien $\h_{j}\cap i
\h_{j}=\h_{j}$ et $\h_{j}$ est une sous-alg\`ebre simple complexe
de $\g$. On remarquera que cette deuxi\`eme possibilit\'e est
exclue, si $B$ est sp\'eciale, puique $\h_{j}$ serait alors
semi-simple complexe et isotrope pour $B$.
\ste On suppose que, pour $j=1,\dots , p$,  $\h_{j}\cap i
\h_{j}=\{0\}$, et que pour $j=p+1,\dots,r$, $\h_{j}\cap i
\h_{j}=\h_{j}$. Si $j=1,\dots,p$,  on note $\k_{j}= \h_{j}\oplus i
\h_{j}$. Si
$j=p+1,\dots,r$,  on note $\k_{j}$ la somme des projections de
$\h_{j}$ dans les id\'eaux simples de $\m$. On note aussi $\k'_{j}=
\h_{j}+i
\h_{j}$, pour $j=1, \dots, r$. On note
$\k=\sum _{j=1, \dots, r}\k_{j}$ et $\k'=\sum_{j=1,\dots,
r}
\k'_{j}= \h+\i\h$ , qui est contenu dans $\k$.  Par ailleurs,
deux \'el\'ements,
$X$ et
$Y$, de
$\m$ commutent si et seulement $X$ et $iY$ commutent (resp. $X$
commute \`a chacune des projections de $Y$ dans les id\'eaux simples de
$\m$). Il en r\'esulte que, pour $j\not=l$,
$\k_{j}$ commute \`a $\h_{l}$, donc $\k_{j}\cap (\sum_{l\not=
j}\k_{l})$ est contenu dans le centre de $\k_{j}$
, qui est semi-simple complexe.   Il en
r\'esulte que :
\begin{equation}\k=\oplus_{j=1, \dots,r}
\k_{j},\>\> \k'=\oplus_{j=1,
\dots,r}\k'_{j}
\end{equation}
Alors $\k$, $\k'$  sont des   sous-alg\`ebres de Lie semi-simples
complexes de $\m$. \ste
Montronsque : \beq\k\cap \n'=\{0\}\eeq
En effet $\n'$ est un id\'eal dans $\h+\n'$, puisque $\i'=\i\cap\m$
et $\n'$ est l'intersection de l'id\'eal $\n$ de $\i$ avec $\m$.
C'est donc un $\h$-module, et aussi un $\k'$-module puisque $\n'$
est un espace vectoriel complexe. Donc
$\k\cap\n'$ est  un sous-$\k'$-module, et aussi une
sous-alg\`ebre r\'esoluble  de $ \k$. Il est clair que les $\k_{j}$
sont des sous-$\k'$-modules de  $\k$, qui n'ont aucun sous-quotient
simple en commun. En effet, d'une part  $\k'_{l}$ agit trivialement
sur
$\k_{j}$, si $j\not=l$. D'autre part, d'apr\`es les d\'efinitions, on
voit que  les sous-quotients simples du
$\k'_{j}$-module $\k_{j}$ sont  sous-quotients  de
$\k'_{j}$, dont aucun n'est trivial, puique $\k'_{j}$ est
une alg\`ebre de Lie semi-simple. Donc, si
$\k\cap
\n'$ est non nul, il a une intersection non nulle,
$\k''$, avec l'un des $\k_{j}$, qui est un $\k'_{j}$-sous-module.
 Comme $\k_{j}$ est une sous-alg\`ebre de Lie de $\m$, il en va de
m\^{e}me de $\k''$, qui est de plus r\'esoluble, puisque c'est le cas de
$\n'$.\ste  Si
$j=1,\dots ,p$,
$\k_{j}=\k'_{j}$ et un $\k'_{j}$-sous-module de $\k_{j}$ est un
id\'eal de $\k_{j}$. Alors $\k\cap
\n'$ est \`a la fois semi-simple et r\'esoluble. Une contradiction
qui montre (1.9) dans ce cas.\ste Si $j=p+1,\dots ,r$,
$\k'_{j}=\h_{j}$ est simple, donc l'une des projections de $\k''$
sur un id\'eal simple de $\m$ est \'egale \`a $\h_{j}$. Cette projection
\'etant un morphisme d'alg\`ebres de Lie, il en r\'esulte que l'alg\`ebre
de Lie r\'esoluble $\k'' $, admet un quotient semi-simple. Une
contradiction qui ach\`eve de prouver (1.9). \ste
Pour $j=p+1, \dots, r$, $\h_{j}$ ne peut \^{e}tre contenu dans un
id\'eal simple de $\m$. En effet, d'apr\`es le Corollaire du Lemme  2,
la restriction de $B$ \`a $\m$ est une forme de Manin. D'apr\`es le
Lemme 1 (ii) et le Lemme 2 (v), la restriction de $B$ \`a un id\'eal
simple de $\m$ est sp\'eciale, et notre assertion en r\'esulte, car
pour $j=p+1, \dots,r$, $\h_{j}$ est isotrope et semi-simple
complexe . Pour $j=p+1, \dots,r$, on notera
$n_{j}$, le nombre d'id\'eaux simples de
$\m$ dans lesquels
$\h_{j}$ a une projection non nulle. On vient de voir que :
\beq n_{j}\geq 2, \>\> j=p+1, \dots, r\eeq
Montrons que
$\n'=\{0\}$.  \ste
On a \'evidemment :
$$dim_{\R}\i'=(\sum_{j=1,\dots,r}dim_{\R}\h_{j})+dim_{\R}\n'$$

Alors,  en posant $p_{j}=1$, pour $ j=1,\dots,p$ et $p_{j}=2$, pour
$ j=p+1,\dots,r$, on a :
\beq dim_{\R}\i'=(\sum_{j=1,\dots,r}p_{j}dim_{\C}\k'_{j})
+dim_{\R}\n'
\eeq
Par ailleurs on a, gr\^{a}ce \`a  (1.7) et (1.9)  :
\beq
dim_{\R}\i'=dim_{\C}\m\geq dim_{\C}\k+dim_{\C}(\n')
\eeq
En posant $n_{j}=1$ pour $j=1,\dots,p$, on a imm\'ediatement :
\beq
dim_{\C}(\k)=\sum_{j=1,\dots,r}n_{j}dim_{\C}\k'_{j}
\eeq
Comme $n_{j}$ est sup\'erieur o\`u \'egal \`a $p_{j}$, pour tout $j$,
d'apr\`es  (1.10), on d\'eduit que l'in\'egalit\'e dans (1.12) est une
\'egalit\'e et que : $$\m=\k\oplus \n'
\>\>et\>\>p_{j}=q_{j},\>\>j=1,\dots,r
$$  La projection de
$\k$ sur un id\'eal simple de $\m$ est \'egal \`a celle de $\k'$,
d'apr\`es la d\'efinition de
$\k$. La projection de l'\'egalit\'e ci-dessus montre que chaque id\'eal
simple de $\m$,  est la somme (pas n\'ecessairement
directe) de la projection de $\k'$ et de celle de $\n'$. Comme on
voit facilement que $\n'$ est un id\'eal de $\k'+\n'$, la
projection de $\n'$  appara\^{i}t comme un id\'eal r\'esoluble de cette
alg\`ebre de Lie semi-simple, donc elle est nulle. Il en r\'esulte que :
$$\n'=\{0\}$$
comme d\'esir\'e. En outre $\m=\k$, donc les $\k_{j}$ sont des id\'eaux
de $\m$. On pose
${\tilde \m_{0}}=\oplus_{j=1,\dots,p}\k_{j}$,
${\tilde \h_{0}}=\oplus_{j=1,\dots,p}\h_{j}$ . On pose $q= r-p$.
Pour
$l=1,
\dots, q$,
$\k_{l}$ est somme de deux id\'eaux simples, $\m'_{l}$, $\m''_{l}$,
car
$n_{p+l}=2$. La projection de $\h_{l+p}$ sur chacun de ces id\'eaux
est bijective, sa surjectivit\'e r\'esultant de la d\'efinition de
$\k_{l}$, son injectivit\'e r\'esultant de la simplicit\'e de
$\h_{l+p}$ et de la non nullit\'e de ce morphisme d'alg\`ebres de
Lie. Donc
$\h_{p+l}:=\{(X,
\tau_{l}(X))\vert X \in \m'_{l}\}$, o\`u $ \tau_{l}$  est un
isomorphisme  $\C$-lin\'eaire de l'alg\`ebre de Lie  $\m'_{l}$ sur
$\m''_{l}$. Donc $\h$  a la forme voulue \qed

\begin {lem}
Aucun poids non nul de $\a$  dans
$\g$ n'est nul sur $\i_{\a}$.\ste

\end{lem}
{\em D\'emonstration :} Raisonnons par l'absurde et supposons
qu'il existe un poids non nul $\alpha$ de $\a$ dans $\g$, nul sur
$\i_{\a}$. Soit
$H_{\alpha}\in \a$ tel que :
\begin{equation}
K_{\g}(H_{\alpha}, X)=\alpha (X),\>\> X\in \a
\end{equation}
Alors
$H_{\alpha}$ appartient \`a l'un des id\'eaux simples de $\g$. En
effet, soit $\j$ une sous-alg\`ebre de Cartan de $\g$, contenant
$\a$. Alors $\alpha $ est la restriction \`a $\a$ d'une racine $\beta$
de $\j$ dans $\g$ et l'on a :
\begin{equation}
K_{\g}( H_{\alpha}, H_{\alpha})>0
\end{equation}
 Soit
$H_{\beta}\in \j$ tel que :
$$K_{\g}(H_{\beta}, X)=\beta (X),\>\> X\in \j$$
Alors $\j$ (resp. $\a$) est la somme directe de ses intersections
avec les id\'eaux simples de $\g$, et
$H_{\beta}$ (resp. $H_{\alpha}$) appartient \`a l'une de celles-ci. On
d\'eduit alors du Lemme 1 (i), qu' il existe $\mu\in \C$ , non nul
car $B$ est non d\'eg\'en\'er\'ee, tel que :
\begin{equation}
B(\lambda H_{\alpha}, X)=Im ( K_{\g}(\lambda \mu
H_{\alpha}, X)),
\>\> \lambda
\in \C \>\> X\in
\g
\end{equation}
Comme $\alpha$ est nulle sur $\i_{\a}$,   il r\'esulte
de (1.14) et (1.16) que $\C H_{\alpha}$ est orthogonale \`a $\i_{a}$,
pour $B$. Comme $B$ est une forme de Manin,  la restriction de $B$ \`a
$\a$ est de signature $(dim_{\C}\a, dim_{\C}\a)$. Tenant compte de
(1.7), on voit que $\i_{\a}$ est un sous-espace de $\a$, isotrope
pour $B$, de dimension maximale. Alors, ce qui pr\'ec\`ede montre que
$\C H_{\alpha}$ est contenu dans $\i_{\a}$. Par ailleurs,  si
$\lambda$ est une racine carr\'ee de $i\mu^{-1}$, $B( \mu H_{\alpha},
\mu H_{\alpha})$ est non nul d'apr\`es (1.15) et (1.16). Une
contradiction avec le fait que $\i_{\a}$ est isotrope qui ach\`eve de
prouver le Lemme. \qed
{\em Fin
de la d\'emonstration du Th\'eor\`eme 1: }\ste
Montrons la propri\'et\'e suivante :
\begin{eqnarray} &&Toute\>\> sous-alg\grave{e}bre\>\>
parabolique\>\>
\q \>\>de\>\>
\g\>\> est\>\> \acute{e}gale \>\>au \>\>normalisateur \nonumber \\&&
dans\>\>
\g \>\>de \>\>son\>\> radical \>\>nilpotent \>\>\w\>\>
\end{eqnarray}   D'abord $\q$ normalise $\w$. Donc,
le normalisateur $\r$ de $\w$ dans $\g$  contient $\q$. C'est
donc une sous-alg\`ebre parabolique de $\g$, qui contient $\w$
comme id\'eal. Compte tenu de [Bou], Ch. I, Paragraphe 5.3, Remarque 2, $\w$
est contenu dans le radical nilpotent, $\x$, de $\r$. On a alors
: $\p\subset \r,\>\>\w\subset \x$, d'o\`u l'on d\'eduit facilement
que $\x=\w$ et $ \r=\q$, comme d\'esir\'e. Donc $\q$ est bien le
normalisateur de $\n$. Ceci ach\`eve de prouver (1.17).\ste
Montrons que $\n$ est le radical nilpotent de $\i$. En effet, comme
$\n$ est somme de sous-espaces poids sous $\a$ et qu'aucun de ces
poids n'est nul sur $\i_{\a}$, d'apr\`es le Lemme pr\'ec\'edent, on a :
$$[\i_{\a}, \n]= \n$$ Donc $[\i,\i]=\h\oplus \n$ et l'intersection
de $[\i,\i] $ avec le radical $\r=\i_{\a}+\n$ de $\i$ est \'egal \`a
$\n$. Donc, d'apr\`es [Bou] Ch. I, Paragraphe 5.3, Th\'eor\`eme 1, $\n$ est
bien le
radical nilpotent de $\i$. On a donc montr\'e que
$\i$ s'\'ecrit de la mani\`ere voulue, pour une d\'ecomposition de
Langlands particuli\`ere  du normalisateur, $\p$, de $\n$ .  Si
$\p=\l'\oplus
\n$ est une autre d\'ecomposition de Langlands de $\p$, $\l$ et $\l'$
sont isomorphes, puisqu'elles sont toutes les deux isomorphes \`a
$\g/\n$. Comme $\i$ contient $\n$, les intersections de $\i$
avec $\l$ et $\l'$ se correspondent dans cet isomorphisme, et la
d\'ecomposition de $\i$ qu'on en d\'eduit, relativement \`a cette
nouvelle d\'ecomposition de Langlands de $\p$, a les propri\'et\'es
voulues.\ste
Etudions la partie r\'eciproque du Th\'eor\`eme. Une sous-alg\`ebre
parabolique de $\g$ est la somme de ses intersections avec
les id\'eaux simples de $\g$. En outre, elle est orthogonale \`a son
radical nilpotent, pour la forme de Killing de $\g$. On conclut
que si $\i$ a une d\'ecomposition comme dans l'\'enonc\'e, elle est
isotrope pour $B$, et de dimension r\'eelle \'egale \`a  la dimension
complexe de $\g$. \qed

\begin{defi}On rappelle qu'une \sac   d'une alg\`ebre
de Lie semi-simple r\'eelle est une \sac fondamentale si et seulement
si elle contient des \'el\'ements r\'eguliers dont l'image par la
repr\'esentation adjointe n'a que des valeurs propres imaginaires
pures. Cela \'equivaut au fait qu'aucune racine de cette \sac
n'est r\'eelle.\ste  Une \sac d'une alg\`ebre de Lie r\'eelle est dite
fondamentale si sa projection dans une sous-alg\`ebre de Levi,
parall\`element au radical, est une \sac fondamentale de cette
alg\`ebre de Lie semi-simple r\'eelle.
\end{defi}
Comme toutes les sous-alg\`ebres de Cartan  fondamentales d'une
alg\`ebre de Lie r\'eelle semi-simple sont conjugu\'ees entre elles par
des automorphismes int\'erieurs, il en va de m\^{e}me pour les
sous-alg\`ebres de Cartan fondamentales d'une alg\`ebre de Lie r\'eelle.
En effet il suffit d'adapter la preuve du fait que
 (ii) implique (i), dans [Bou], Ch. VII, Paragraphe 3.5,
 Proposition 5, en
remarquant pour cela que tout automorphisme int\'erieur d'une alg\`ebre
de Levi d'une alg\`ebre de Lie r\'eelle, s'\'etend en un automorphisme
int\'erieur de l'alg\`ebre de Lie.

\begin{lem}
On conserve les hypoth\`eses et notations du Th\'eor\`eme 1. \ste
(i) Si $\f$ est une sous-alg\`ebre de Cartan de $\h$, il existe
des \'el\'ements r\'eguliers de $\g$ contenus dans $\f\oplus
\i_{\a}$. Le centralisateur dans $\g$, $\j$, de ${\tilde
\f}:=\f\oplus
\i_{\a}$
est une sous-alg\`ebre de Cartan
de $\g$, contenue dans $\l$, v\'erifiant $\j=(\j\cap\m)\oplus
\a$.\ste (ii) Si
$\f$ est une sous-alg\`ebre de Cartan de $\h$,
$ \f\oplus
\i_{\a}$ est une sous-alg\`ebre de Cartan  de $\i$.
\ste (iii) Toute
sous-alg\`ebre de Cartan de $\i$ (resp. sous-alg\`ebre de
Cartan de $\i$ contenue dans $\h+\i_{\a}$) est conjugu\'ee, par un
automorphisme int\'erieur de
$\i$,  (resp. \'egale)   \`a une alg\`ebre de ce type.\ste
(iv) Soit ${\tilde \f'}$ une \sac de $\i$. Il existe une unique
d\'ecomposi-\ste tion de Langlands $\l'+\n$ de $\p$, telle que $\l'$
contienne ${\tilde \f'}$. Alors, notant  $\f':={\tilde
\f'}\cap\l'^{der}$, on a ${\tilde \f'}=\f'+(\i\cap\a')$, o\`u $\a'$
est le centre de $\l'$. De plus ${\tilde \f'}$ est une \sac
fondamentale de $\i$, si et seulement si $\f'$ est une \sac
fondamentale de
$\h':=\i\cap\l'^{der}$
\end{lem} \dem
Montrons (i). On raisonne par l'absurde. On note $\j'$ une \sac de
$\g$, qui contient $\f\oplus \a$. Celle-ci  existe puisque les
\'el\'ements de  $\f\oplus \a$ sont semi-simples. Supposons qu'aucun
\'el\'ement  de
${\tilde \f}:=\f\oplus
\i_{\a}$ ne soit r\'egulier dans
$\g$. Alors, pour tout $X\in {\tilde \f}$, il existe une racine
$\alpha_{X} $ de $\j'$ dans $\g$, nulle sur $X$. Pour une racine
donn\'ee, l'intersection de son noyau avec ${\tilde \f}$ est un ferm\'e
de ${\tilde \f}$. Notre hypoth\`ese montre que ${\tilde \f}$ est la
r\'eunion de ces ferm\'es. Il en r\'esulte que l'un de ces sous-espaces
vectoriels est  d'int\'erieur non vide, donc \'egal \`a ${\tilde \f}$.
Cela signifie qu'une racine de $\j'$ s'annule sur ${\tilde \f}$.
Alors, d'apr\`es le Lemme 7, dont on v\'erifie ais\'ement qu'il  est
valable pour toute d\'ecomposition de Langlands de $\p$ , celle-ci
doit \^{e}tre nulle sur
$\a$. C'est donc une racine de
$\j'$ dans $\m$, qui ne peut \^{e}tre nulle sur $\f$. Une contradiction
qui prouve la premi\`ere partie de (i). Le centralisateur $\j$ de
${\tilde \f}$ est donc une sous-alg\`ebre de Cartan. Par ailleurs
$\j$ contient  $\a$. On en d\'eduit la deuxi\`eme partie de (i).
\ste \ D'apr\`es le Lemme 7,  le nilespace de ${\tilde \f}$ dans $\n$
est r\'eduit \`a z\'ero. Comme $\f$ est une sous-alg\`ebre de
Cartan de $\h$, le nilespace de ${\tilde \f}$ dans $\h$ est
\'egal \`a $\f$. Finalement le nilespace de ${\tilde \f}$ dans
$\i$ est \'egal \`a ${\tilde \f}$. Alors (ii) r\'esulte de [Bou], Ch. VII,
Paragraphe 2.1, Proposition 3. \ste Montrons (iii). Soit ${\tilde
\f}'$ une
autre sous-alg\`ebre de Cartan de $\i$. La projection, ${\tilde
\f''}$, de
${\tilde
\f}'$ sur
${\tilde \h}:=\h+\i_{\a}$, parall\`element \`a $\n$, est une
sous-alg\`ebre de Cartan de ${\tilde \h}$ (cf. [Bou], Ch. VII,
 Paragraphe 2.1,
Corollaire 2 de la Proposition 4), donc de la forme $\f'+\i_{\a}$,
o\`u $\f'$ est une sous-alg\`ebre de
Cartan de $\h$. Alors ${\tilde \f}'$ et ${\tilde \f}''$ sont deux
sous alg\`ebres de Cartan de $\i$, ayant la m\^{e}me projection sur
${\tilde \h}$, parall\`element \`a $\n$, donc conjugu\'ees par un
automorphisme int\'erieur de $\i$, d'apr\`es [Bou], Ch. VII,
Paragraphe 3.5,  Proposition 5  (voir aussi
apr\`es la D\'efinition 3). Si de plus ${\tilde \f}'$ est contenue dans
$ {\tilde
\h}$, le raisonnement ci-dessus montre qu'elle a la forme indiqu\'ee.
 Ce
qui prouve (iii). \ste
Prouvons (iv). Gr\^{a}ce \`a (iii), on se ram\`ene, par conjugaison,
au  cas o\`u
${\tilde
\f'}$ est contenue dans $\l$ et comme dans (i). Si $\l'+\n$ est une
d\'ecomposition de Langlands  de $\p$, o\`u $\l'$ contient ${\tilde
\f'}$, $\l'$ contient un
\'el\'ement r\'egulier de $\g$, contenu dans ${\tilde
\f'}$, dont le centralisateur dans $\g$ est une \sac de $\g$,
contenue dans $\l'$. Celle-ci est \'egale au centralisateur dans $\g$
de ${\tilde
\f'}$. D'o\`u l'unicit\'e de $\l'$, gr\^{a}ce aux propri\'et\'es des
d\'ecompositions de Langlands (cf. Lemme 3 (i)). L'assertion
sur les \sac fondamentales est claire car $\h'$ est une
sous-alg\`ebre de Levi de $\i$, d'apr\`es le Th\'eor\`eme 1.\qed

\section{Triples de Manin pour une alg\`ebre de Lie r\'eductive
complexe : Descente
 }
\setcounter{equation}{0}
Dans toute la suite $\g$ d\'esignera une alg\`ebre de Lie
r\'eductive complexe. On fixe,  $\j_{0}$, une sous alg\`ebre de
Cartan de
$\g$, $\b_{0}$ une sous-alg\`ebre de  Borel de $\g$, contenant
$\j_{0}$. On note
$\b'_{0}$ la sous-alg\`ebre de Borel oppos\'ee \`a $\b_{0}$,
relativement \`a $\j_{0}$.
\begin{defi} Un triple de Manin pour $\g$ est un triplet
$(B,
\i,\i')$, o\`u $B$ est une forme de Manin  sur $\g$, $\i$ et $\i'$
sont des sous-alg\`ebres de Lie r\'eelles de $\g$,  isotropes pour $B$,
telles que $\g=:\i\oplus \i' $. La signature de $B$ \'etant \'egale \`a
$(dim_{\C}\g,dim_{\C}\g) $,  $\i$ et
$\i'$ sont Lagrangiennes. Un s-triple  est un triple de Manin
pour lequel la forme est sp\'eciale. \ste  Si $\i$ est sous $\p$ et
$\i'$ est sous
$\p'$, on dit que le triple de Manin est sous $(\p,\p')$

\end{defi}
\begin{rem}   D'apr\`es le Lemme 2 (v), si $\g$ est simple, la notion
de s-triple  et de triple de Manin coincident.
\end{rem}
On note $G$ le groupe  connexe , simplement
connexe, d'alg\`ebre de Lie $\g$. Si $\s$ est une sous-alg\`ebre de
$\g$, on note $S$ le sous-groupe analytique de $G$, d'alg\`ebre de
Lie $\s$. Comme $\g$ est complexe, les sous-groupes
paraboliques de $\g$ sont connexes (cf. [Bor], Th\'eor\'eme
11.16). Donc, si $\p$ est une sous-alg\`ebre parabolique de $\g$,
$P$ est le sous-groupe parabolique de $G$, d'alg\`ebre de Lie
$\p$.  \ste On remarque  que $G$ agit sur l'ensemble des
triples de Manin , en posant, pour tout triple de Manin
$(B,\i,\i')$ et tout
$g\in G$ : $$g (B,\i,\i'):=(B,Ad\>\>g(\i),Ad\>\>g(\i'))$$
  Notre  but est construire, par r\'ecurrence sur la dimension de
$\g^{der}$, tous les triples de Manin modulo cette action de $\g$.
\begin{prop}
Tout triple de Manin est conjugu\'e, sous l'action de $G$, \`a un
triple de Manin
 $(B,\i,\i')$  sous $(\p,\p')$ , avec $\b_{0}\subset \p$ et
$\b'_{0}\subset
\p'$ (un tel triple de Manin
 sera dit standard ).\ste
De plus $\p$ et $\p'$ sont uniques.
\end{prop}
\dem
Montrons d'abord  que :$$ \>\>L\>\>'\>\>intersection \>\>de\>\>
deux \>\>sous-alg\grave{e}bres\>\> de \>\>Borel\>\> de\>\> \g\>\>,
\>\>{\underline
\b}, \>\>{\underline \b'},$$\beq \>\> contient
\>\> une \>\>sous-alg\grave{e}bre \>\>de\>\> Cartan \>\>de\>\> \g
\eeq\ste
On a d'abord ${\underline \b'}= Ad\> g( {\underline \b})$, pour
un \'el\'ement $g$ de $G$. Soit
${\underline \j}$  une sous-alg\`ebre de Cartan de $\g$, contenue
dans
${\underline \b}$ et
$W$ un ensemble de repr\'esentants dans $\g$ du groupe de Weyl de la
paire
$(\g,{\underline \j})$. La d\'ecomposition de Bruhat implique qu'il
existe
$b$,
$b_{1}\in {\underline B}$ et
$w\in W$ tels que :
$g=bwb_{1}$. Alors $\j:=Ad\>b ({\underline \j})$ est contenue dans
${\underline \b}$, puisque ${\underline B}$ normalise ${\underline
\b}$. Comme $W$ normalise
${\underline
\j}$, on a aussi $\j= Ad\>b \>(w({\underline \j}))$. Mais
${\underline \b'}=Ad\>bw \>({\underline \b})$, puisque $b_{1}\in
{\underline B}$, normalise ${\underline \b}$. Finalement
$\j$ est contenue dans ${\underline \b}\cap {\underline \b'}$,
comme d\'esir\'e.\ste Soit $(B, {\underline \i}, {\underline \i'}) $
 un
triple de Manin sous
$({\underline \p}, {\underline \p'})$. Soit ${\underline \b}$
(resp. ${\underline \b'}$) une sous-alg\`ebre de Borel de $\g$,
contenue dans
${\underline \p}$ (resp. ${\underline \p'})$. \ste
On a $\g={\underline \i}+{\underline \i'}\subset {\underline
\p}+ {\underline \p'}$. Donc ${\underline
\p}+ {\underline \p'}$ est \'egal \`a $\g$ et ${\underline
P}\>{\underline P'}$ est ouvert dans $G$. Mais ${\underline
P}\>{\underline P'}$ est r\'eunion de $({\underline B}, {\underline
B'})$-doubles classes, qui sont en nombre fini (Bruhat). L'une de
ces doubles classes contenues dans ${\underline
P}\>{\underline
P'}$ doit donc \^{e}tre ouverte. Soit $p\in {\underline
P}$ et $p' \in {\underline
P'}$, tels que ${\underline B}pp'{\underline B'}$ soit un ouvert  de
$G$. On pose
$B_{1}=p^{-1}{\underline B}p$, $B'_{1}=p'{\underline B'}p'^{-1}$.
Alors le sous-groupe de Borel de $G$, $B_{1}$ (resp.
$B'_{1}$), est contenu dans $P$ (resp. $P'$) et $B_{1}B'_{1}$ est
ouvert dans $G$. Donc, on a $\b_{1}+\b'_{1}=\g$ et l'intersection
de
$\b_{1}$ et $ \b'_{1}$ contient une  sous-alg\`ebre de Cartan de
$\g$, $\j_{1}$ (voir ci-dessus). Pour des raisons de dimension,
cette intersection est r\'eduite \`a $\j_{1}$. Alors $\b_{1}$ et $
\b'_{1}$ sont oppos\'ees relativement \`a $\j_{1}$ . D'apr\`es [Bor],
Proposition 11.19, il existe $g'\in G$ tel que
$Ad\>g'(\b_{1})=\b_{0}$, $Ad\>g'(\j_{1})=\j_{0}$. Alors
 $Ad\>g'(\b'_{1})$ est \'egal \`a $\b'_{0}$. Alors, notant
 $\i=Ad\> g'(  {\underline \i})$, $\i'=Ad\> g'(  {\underline \i'})$,
on voit que $(B,\i,\i') $ v\'erifie les propri\'et\'es voulues.\ste
L'unicit\'e de $\p$ r\'esulte du fait que deux sous-alg\`ebres
paraboliques de $\g$, conjugu\'ees par un \'el\'ement de $G$ et contenant
une m\^{e}me sous-alg\`ebre de Borel, sont \'egales (cf. [Bor],
Corollaire 11.17).\qed
\ste
{\em On fixe d\'esormais $\p$ (resp. $\p'$) une sous-alg\`ebre
parabolique de
$\g$, contenant $\b_{0}$ (resp. $\b'_{0})$. On note $\p=\l\oplus \n$
(resp. $\p'=\l'\oplus \n'$) la d\'ecomposition de Langlands de $\p$
(resp. $\p'$) telle que $\l$ (resp. $\l'$) contienne $\j_{0}$ (cf.
Lemme 3).  On note $\m=\l^{der}$, $\a$ le centre de $\l$. Si $\i$
est une sous-alg\`ebre de Lie r\'eelle de $\g$, Lagrangienne  pour une
forme de Manin, on notera $\h=\i\cap \m$, $\i_{\a}=\i\cap\a$,
${\tilde
\h}=\h\oplus\i_{a}$. On introduit des notations similaires
pour $\p'$.}\ste  Comme $\b_{0}\subset \p$
(resp. $\b'_{0}\subset \p' $), $\n$ (resp. $\n'$) est contenu dans
le radical nilpotent de $\b_{0}$ (resp. $\b'_{0}$). Ces derniers
sont d'intersection r\'eduite \`a z\'ero, donc :
\begin{equation}
\n\cap\n'=\{ 0\}
\end{equation}
D\'ecomposant $\p\cap\p'$ en sous-espaces poids sous $\j_{0}$, on
voit que :
\begin{equation}
\p\cap\p'=(\l\cap\l')\oplus (\n\cap\l')\oplus (\n'\cap\l).
\end{equation}

\begin{prop}
(i) Si un \'el\'ement de $G$ conjugue deux triples de Manin  sous $(\p,
\p')$, c'est un \'el\'ement de $P\cap P'$.\ste
(ii) Le groupe $L\cap L'$ est \'egal au sous-groupe analytique de
$G$, d'alg\`ebre de Lie $\l\cap \l'$. Notons $N_{L'}$ (resp. $N'_{L}$)
, le sous-groupe analytique de $G$, d'alg\`ebre de Lie $\n\cap\l'$
(resp. $\n'\cap\l$). Alors on a :
$$P\cap P'=( L\cap L')N_{L'}N'_{L}$$
De plus $N_{L'}$ et $N'_{L}$ commutent entre eux.
\end{prop}
\dem  Si $(B,\i,\i')$ et $(B, {\underline \i}, {\underline \i'})
$ sont deux triples de Manin sous $(\p, \p')$, conjugu\'es par un
\'el\'ement,
$g$, de
$G$, celui-ci conjugue le radical nilpotent de $\i$ avec celui de
${\underline  \i}$, donc normalise $\n$, puisque les deux triples
de Manin sont sous $(\p,\p')$. Mais un \'el\'ement du normalisateur,
$Q$, dans
$G$ de
$\n$, normalise le normalisateur dans $\g$ de $\n$, c'est \`a dire
$\p$, comme on l'a vu  plus haut ( cf.  (1. 17)).
Comme $P$ est connexe, les \'el\'ements de $Q$ normalisent $P$. Donc
$Q$ est inclus dans $P$ et $g\in P$. De m\^{e}me, on a  $g\in P'$. D'o\`u
(i)\ste
Montrons (ii). Il est clair que $P\cap P'$ est un sous-groupe de Lie
de $G$, d'alg\`ebre de Lie $\p\cap\p'$.
On a : $$[\n\cap\l',\n'\cap\l]\subset [\n,\l]\cap
[\n',\l']\subset\n\cap \n'$$
Donc $N_{L'}$ et $N'_{L}$ commutent entre eux, d'apr\`es (2.1). Alors
$ ( L\cap L')^{0}N_{L'}N'_{L}$ est un sous-goupe ouvert et connexe
de $P\cap P'$, donc on a :
\begin{equation}
(P\cap P')^{0}=( L\cap L')^{0}N_{L'}N'_{L}
\end {equation}
Soit $g\in P\cap P'$. Alors $Ad\>g (\j_{0})$ est une sous-alg\`ebre de
Cartan de $\g$, contenue dans $\p\cap\p'$, c'est donc une
sous-alg\`ebre de Cartan de
$\p\cap\p'$ (cf [Bou], Ch. VIII, Paragraphe 2.1, Exemple 3), donc conjugu\'e,
par un \'el\'ement $g'$ de $(P\cap P')^{0}$, \`a $\j_{0}$, puisqu'il
s'agit d'alg\`ebres de Lie complexes. Donc
$ Ad\>g' g(\j_{0})=\j_{0}$ et $g'g$ est un \'el\'ement de $P\cap P'$. En
utilisant la d\'ecomposition de Bruhat de $G$ et $P$, pour $B_{0}$,
on voit que $g'g$ centralise le centre $\a$ de $\l$. De m\^{e}me
on voit que $g'g$ centralise le centre $\a'$ de $\l'$. Donc $g'g$
est un \'el\'ement du centralisateur, $L''$,  de $\a+\a'$ dans $G$. Mais
$P'':=L''(N'_{L} N_{L'})$, est une d\'ecomposition de Langlands du
sous-groupe parabolique de $G$, d'alg\`ebre de Lie :
$$\p''=(\p'\cap\l)\oplus \n= (\l\cap\l')\oplus (\l\cap\n')\oplus
\n$$ Or $P''$ est connexe, puisque $G$ est complexe. Donc $L''$ est
connexe. Par ailleurs, il  contient $L\cap L'$ et a m\^{e}me alg\`ebre
de Lie que $L\cap L'$. Donc on a :
\begin{equation}
 L''=L\cap L'= (L\cap L')^{0}
\end{equation}
On conclut alors que $g'g\in (L\cap L')^{0}$. Donc $g$ est un
\'el\'ement de $ (P\cap P')^{0}$. Ce qui pr\'ec\`ede montre que :

$$P\cap P'=(P\cap P')^{0}$$
On ach\`eve de prouver $(ii)$, gr\^{a}ce \`a (2.4) et (2.5).\qed
Le Lemme suivant est une cons\'equence facile de r\'esultats de
Gantmacher (cf. [G]).
\begin{lem}
Si $\sigma$ et $\sigma'$ sont deux automorphismes involutifs et
antilin\'eaires d'une alg\`ebre de Lie semi-simple complexe, $\m$,
celle-ci contient au moins un \'el\'ement non nul et  invariant par ces
deux involutions
\end{lem}
\dem
Avec nos hypoth\`eses
$\sigma\sigma'$ est un automorphisme
$\C$-lin\'eaire de
$\m$,
dont l'espace des points fixes, $\m^{
\sigma\sigma '}$, est un espace vectoriel complexe,  non r\'eduit \`a
z\'ero d'apr\`es  [G], Th\'eor\`eme  28. Mais $ \m^{
\sigma\sigma '}$ , est \'egal \`a $\{ X\in \m\vert
\sigma(X)=\sigma' (X)\}$ donc aussi \'egal \`a $\m^{
\sigma '\sigma }$. Si $X \in \m^{
\sigma\sigma '}$, on a donc $\sigma '(\sigma(X))= X$, soit encore
$\sigma '(\sigma(X))=\sigma (\sigma(X))$. Donc $\sigma(X)$ est
\'el\'ement de $\m^{
\sigma\sigma '}$. Par suite $\sigma$, restreint \`a $\m^{
\sigma\sigma '}$ est une involution antilin\'eaire de $\m^{
\sigma\sigma '}$. L'ensemble de ses points fixes est  une forme
r\'eelle de $\m^{
\sigma\sigma '}$, donc il est non r\'eduit \`a z\'ero. Mais cet ensemble
est \'egal \`a $\m^{\sigma}\cap\m^{\sigma '}$. \qed

\begin{prop}
Si $\sigma$ et $\sigma'$ sont deux af-involutions d'une alg\`ebre
de Lie semi-simple complexe, $\m$, elle contient au moins un \'el\'ement
non nul et  invariant par ces deux involutions.
\end{prop}
\dem
On note
$\m_{j}$,
$j= 1,
\dots,r$, les id\'eaux simples de $\m$. On d\'efinit une involution
$\theta$ de
$\{1,
\dots,r\}$ caract\'eris\'ee par  : $\sigma (\m_{j})=\m_{\theta (j)},
\>\> j =1,
\dots,r$. Nous allons d'abord \'etudier le cas suivant :
\begin{equation}
Il\>\> existe \>\> j \>\>tel\>\> que \>\>\theta(j)=\theta' (j)=j
\end{equation}
Dans ce cas, la restriction de
$\sigma$ et $\sigma'$ \`a $\m_{j}$, sont deux automorphismes
involutifs et
antilin\'eaires de $\m$, d'apr\`es le Corollaire du Lemme 6, qui
ont des points fixes non nuls en commun, d'apr\`es le Lemme
pr\'ec\'edent. La Proposition en r\'esulte, dans ce cas.
\ste Supposons maintenant :
\begin{equation}
Il\>\> existe \>\> j \>\>tel\>\> que \>\>\theta(j)=\theta' (j)
\not=j
\end{equation}
On note $j':=\theta(j)$. Il est clair que :
$(\m_{j}\times \m_{j'})^{\sigma}=\{(X,\sigma(X))\vert X\in \m_{j}\}$
et de m\^{e}me pour $(\m_{j}\times \m_{j'})^{\sigma'}$. Il existe un
\'el\'ement non nul, $X$ de $\m_{j}$ tel que
$({\sigma'}^{-1}{\sigma})(X)=X$, car ${\sigma'}^{-1}{\sigma}$ est
un automorphisme $\R$-lin\'eaire de $\m_{j}$ (cf. [G], Th\'eor\`eme 28).
Alors
$(X,\sigma(X))$ est un  \'el\'ement non nul de $(\m_{j}\times
\m_{j'})^{\sigma}\cap(\m_{j}\times
\m_{j'})^{\sigma'}$ , ce qui prouve la Proposition dans ce cas.
Il nous reste \`a \'etudier le cas suivant :
\begin{equation}
Pour \>\> tout  \>\> j,  \>\>\theta(j)\not=\theta' (j)
\end{equation}
On construit, pour tout $j$, par r\'ecurrence sur $n$, une suite
$j_{1}=j, j_{2},\dots,\ste j_{n},\dots$ , telle que :
\beq
pour\>\> tout\>\> n,\>\> j_{n+1}\not= j_{n}
\eeq
\beq
pour \>\>tout\>\> n,\>\> j_{n+1}= \theta( j_{n}) \>\>ou
\>\>\theta '(
\j_{n})
\eeq
Plus pr\'ecis\'ement, on pose
\beq j_{2}=\theta(1) \>\>si\>\> \theta(1)\not=1,\>\>
j_{2}=\theta'(1) sinon\eeq et, pour
$n\geq 2$, on pose :
$$ j_{n+1}= \theta( j_{n})\>\> si\>\> j_{n}= \theta'(j_{n-1})
\>\>et\>\>\theta (j_{n})\not=j_{n}
$$
$$j_{n+1}= \theta'( j_{n})\>\> si\>\> j_{n}= \theta'(j_{n-1})
\>\>et\>\>\theta (j_{n})=j_{n}$$
\beq
j_{n+1}= \theta'( j_{n})\>\> si\>\> j_{n}= \theta(j_{n-1})
et\>\>\theta' (j_{n})\not=j_{n}
\eeq
$$j_{n+1}= \theta( j_{n})\>\> si\>\> j_{n}= \theta(j_{n-1})
\>\>et\>\>\theta '(j_{n})=j_{n}$$
Ces relations d\'efinissent la suite $(j_{n})$, car,  \`a cause de
(2.8),  on a n\'ecessaire-\ste ment
$\theta(j_{n})\not=\theta'(j_{n})$ et
$\theta(j_{n-1})\not=\theta'(j_{n-1})$. Par ailleurs les relations
(2.9) et (2.10) sont v\'erifi\'ees, la premi\`ere r\'esultant d'une
r\'ecurrence imm\'ediate. On obtient \'egalement les relations suivantes :
$$
Pour\>\> n\geq 2,  si \>\>\theta (j_{n})\not =j_{n}\>\> et
\>\>si\>\>
\theta '(j_{n})\not =j_{n}
\>\>on \>\>a\>\> :$$
\beq( j_{n-1},j_{n}, j_{n+1})\>\>est\>\> \acute{e} gal\>\>
\grave{a}
\>\>(\theta(j_{n}), j_{n},
\theta'(j_{n}))\>\>ou\>\>\grave{a}\>\>
(\theta'(j_{n}), j_{n},
\theta(j_{n}))
\eeq

\beq
Pour \>\> n\geq 2\>\> et \>\> si \>\>\theta (j_{n})=j_{n}\>\> ou
\>\>si\>\>
\theta '(j_{n})=j_{n}
\>\>on \>\> a \>\> : j_{n-1}=j_{n+1}
\eeq
Apr\`es renum\'erotation des $\g_{j}$, on peut supposer que le d\'ebut
de la suite $(j_{n})$, s'\'ecrit $j_{1}=1,j_{2}=2,\dots,
j_{p}=p,j_{p+1}=k<p$
\ste On fait d'abord la convention  suivante :
$$S\>'\> il\>\> existe \>\>j \>\>tel \>\>que \>\>\theta(j)=j
\>\>ou\>\>
\theta'(j)=j,\>\>on \>\>suppose\>\> qu\>'\>on
$$
$$ \>\>l\>'\>a \>\>choisi \>\>comme\>\>
premier\>\> \acute{e}l\acute{e} ment,\>\> et,\>\>
quitte\>\>\grave{a}
\>\> \acute{e} changer
\>\>le\>\> role $$ \beq \>\>de\>\>
\theta\>\> et\>\>
\theta',\>\> qu'il\>\> est\>\> fix \acute{e}\>\> par\>\>
\theta. \>\> On \>\> a \>alors \>\>\theta(1)= 1,\>\>
\theta'(1)=2\eeq
Traitons le cas o\`u $p=2$. Alors $j_{1}=j_{3}$, et
(2.8), (2.13) montrent que $\theta(2)$ ou $\theta'(2)$ est \'egal \`a 2.
Alors on doit avoir
$\theta(1)=1$, d'apr\`es (2.15),  puis $\theta'(1)=2$ d'apr\`es
(2.11). Comme
$\theta'(1)=2$, on a n\'ecessairement
$\theta(2)=2$.  Dans ce cas, un \'el\'ement
$(X_{1},X_{2})\in \m_{1}\oplus\m_{2}$ est invariant par
$\sigma$ et $\sigma'$ si et seulement si on a :
$$X_{1} =\sigma (X_{1} ) ,\>\>       X_{2}= \sigma
(X_{2}),\>\>
X_{2}=\sigma' ( X_{1})$$
ce qui \'equivaut au syst\`eme :
$$X_{1} =\sigma (X_{1} ),\>\> X_{1} =(\sigma'^{-1}\sigma\sigma
')(X_{1} ),\>\>X_{2}=\sigma' ( X_{1})$$
Mais la restriction de $\sigma $  \`a $\m_{1}$
(resp. $\m_{2}$) est un  automorphisme involutif antilin\'eaire,
puisque
$\theta(1)=1$ et $\theta(2)=2$ (cf. le Corollaire du Lemme 6). De
plus, la restriction de
$\sigma' $  \`a $\m_{1}$ est soit
$\C$-lin\'eaire, soit antilin\'eaire, d'apr\`es le Lemme 6. Alors, la
restriction de $\sigma'^{-1}\sigma\sigma
' $  \`a $\m_{1}$ est un  automorphisme involutif antilin\'eaire.
Alors, dans le cas $p=2$, la Proposition r\'esulte du Lemme 12. \ste
On suppose maintenant :
\beq
p>2
\eeq
On remarque d'abord que :
\beq
Si\>\>j=2,\dots, p-1,\>\> on \>\> a \>\> \theta(j)\not=j \>\> et
\>\>\theta'
(j)\not=j
\eeq
En effet, si on avait par exemple $\theta(j)=j$,  (2.14) conduirait
\`a $j-1=j+1$ une contradiction qui prouve (2.17).\ste
Montrons maintenant que : \beq
k=1 \>\> ou \>\> p-1
\eeq
Supposons  $k\not= 1$. Alors, on a $1<k\leq p-1$. Alors d'apr\`es
(2.13) et (2.14), on a l'\'egalit\'e d'ensembles  :

\beq \{\theta(k),\theta'(k)\}=\{k-1, k+1\},
\eeq ce qui  implique :
\beq
\theta(k)\not =k, \>\> \theta'(k)\not =k
\eeq
Comme $j_{p+1}=k$, on d\'eduit de (2.20) et (2.13) que la s\'equence
$(j_{p}, j_{p+1},j_{p+2})$ est \'egale soit \`a
$(\theta(k),k,\theta'(k))$,  soit \`a  $(\theta'(k),k,\theta(k))$,
c'est \`a dire , gr\^{a}ce \`a (2.19)), soit \`a $(k-1, k, k+1)$, soit \`a
$(k+1, k, k-1)$. Mais $j_{p}=p$,
 est diff\'erent de $k-1$. Donc $p=k+1$, i.e. $k=p-1$. Ceci
ach\`eve de prouver (2.18).\ste
Traitons d'abord le cas :
\beq
k=1
\eeq
Comme $p>2$, on a $1\not= p-1$. Donc $j_{p-1}=p-1$, est diff\'erent
de $j_{p+1}=1$. (2.14), (2.13) impliquent l'\'egalit\'e d'ensembles :
\beq \{\theta(p),\theta'(p)\}=\{p-1, 1\}
\eeq
Supposons, d'abord que $\theta'(1)=2$, ce qui implique,
d'apr\`es (2.11), que $\theta(1)=1$. Comme $p>2$, ni $\theta(p)$, ni
$\theta'(p)$ ne peut \^{e}tre \'egal \`a 1. Une contradiction avec
l'\'equation pr\'ec\'edente qui montre que l'on doit avoir, d'apr\`es
(2.11) :
\beq \theta(1)=2\eeq
Comme $p>2$, la seule possibit\'e laiss\'ee par
(2.22) est :
\beq\theta (p-1)=p \>\>et\>\> \theta ' (p)=1\eeq
On d\'eduit de (2.23) et (2.17), joints \`a (2.13), que, pour
$j=1,\dots,p-1$, on a :
\beq \theta(j)=j+1, \>\>si \>\>j\>\> est\>\> impair \>\>(resp.\>\>
\theta'(j)=j+1\>\> si\>\> j\>\> est\>\> pair)
\eeq ce qui, joint \`a (2.24), implique que $p$ est pair. Notons
$p=2q$. \ste
On d\'eduit de (2.25) et (2.23) qu'un \'el\'ement $(X_{1}, \dots,X_{p})$
de
$\m_{1}\oplus\dots\oplus
\m_{p}$, est invariant \`a la fois par $\sigma$ et $\sigma'$ si et
seulement si le syst\`eme suivant est v\'erifi\'e:
$$\sigma(X_{1})=X_{2},\>\>
 \sigma'(X_{2})=X_{3}$$
$$\dots,\>\> \dots $$
$$\sigma(X_{2j-1})=X_{2j},\>\>
 \sigma'(X_{2j})=X_{2j+1}$$
$$\dots, \>\>\dots $$
$$\sigma(X_{2q-1})=X_{2q},\>\>
 \sigma'(X_{2q})=X_{1}$$
Notons $\tau $ la restriction de $(\sigma'\sigma)^{q}$ \`a
$\m_{1}$, qui est un automorphisme $\R$-lin\'eaire de $\m_{1}$.  Ce
syst\`eme poss\`ede une solution non nulle si et seulement si
l'\'equa-\ste tion :
$$X_{1}=\tau (X_{1}), \>\> X_{1}\in \m_{1}$$
poss\`ede une solution non nulle. C'est le cas, d'apr\`es [G],
Th\'eor\`eme 28. Ceci ach\`eve de prouver la Proposition dans le cas
$k=1$.\ste On suppose maintenant :
\beq k=p-1>1\eeq
Comme $(j_{p-1,}j_{p},j_{p+1})=(p-1,p,p-1)$, on d\'eduit
de (2.13), (2.14) et (2.8),  que l'on a  soit :
\beq \theta(p)=p,\>\>\theta'(p)=p-1
\eeq
soit :
\beq \theta(p)=p-1,\>\>\theta'(p)=p
\eeq
Alors, d'apr\`es notre convention  (2.15), on a $\theta(1)=1$.
Supposons (2.27) v\'erifi\'e. Comme ci-dessus, ceci joint \`a (2.23) et
(2.13), montre que $p$ est pair et que , pour $j=1,\dots,p-1$, on a
:
\beq \theta'(j)=j+1, \>\>si \>\>j\>\> est\>\> impair \>\>(resp.\>\>
\theta(j)=j+1\>\> si\>\> j\>\> est\>\> pair)
\eeq
On note $p=2q$.
On d\'eduit de (2.28) et (2.29) qu'un \'el\'ement $(X_{1}, \dots,X_{p})$
de
$\m_{1}\oplus\dots\oplus
\m_{p}$, est invariant \`a la fois par $\sigma$ et $\sigma'$ si et
seulement si le syst\`eme suivant est v\'erifi\'e :
$$\sigma(X_{1})=X_{1},\>\>
 \sigma'(X_{1})=X_{2}$$
$$\dots,\>\> \dots $$
$$\sigma(X_{2j})=X_{2j+1},\>\>
 \sigma'(X_{2j+1})=X_{2j+2}$$
$$\dots, \>\>\dots $$
$$\sigma(X_{2q-2})=X_{2q-1},\>\> \sigma'(X_{2q-1})=X_{2q}$$
 $$\sigma(X_{2q})=X_{2q}$$
Notant $\tau $ la restriction de $(\sigma'\sigma)^{q}$ \`a
$\m_{1}$, qui est un automorphisme $\R$-lin\'eai-\ste re de $\m_{1}$,
ce syst\`eme poss\`ede une solution non nulle si et seulement si
le sys-\ste t\`eme :
\beq X_{1}=\sigma (X_{1}),\>\>   X_{1}=(\tau^{-1}\sigma\tau)
(X_{1}),
\>\> X_{1}\in
\m_{1}\eeq  poss\`ede une solution non nulle. La restriction de
$\sigma$ \`a $\m_{1}$ et $\m_{p}$ est antilin\'eaire. Par ailleurs
$\tau$ est soit $\C$-lin\'eaire, soit antilin\'eaire, d'apr\`es le
Lemme 6.  Donc la restriction \`a $\m_{1}$ de  $\tau^{-1}\sigma\tau$
est antilin\'eaire. Il r\'esulte alors du Lemme 12, que (2.30) \`a une
solution non nulle. Ce qui ach\`eve la preuve de la Proposition dans
le cas \'etudi\'e. Le cas o\`u (2.28) est satisfait se traite de mani\`ere
similaire, mais alors $p$ est impair.\ste
Ceci ach\`eve notre discussion et la preuve de la Proposition.\qed
\begin{theo}
Si $\g$ n'est pas commutative et si  $(B,\i,\i')$ est un triple de
Manin de $\g$, sous $(\p,\p')$,
$\l\cap\l'$ est diff\'erent de $\g$.
\end{theo}
\dem
Raisonnons par l'absurde et supposons qu'il existe un triple de
Manin , $(B,\i,\i')$, sous $(\g,\g)$, et que $\g$ ne soit pas
commutative. Alors $\h:=\i\cap{\g}^{der}$ (resp.
$\h':=\i'\cap\g^{der}$) est  l'espace des points fixes
d'une af-involutions $\sigma$ (resp. $\sigma '$) de
$\g^{der}$, d'apr\`es le Th\'eor\`eme 1. En appliquant la Proposition
pr\'ec\'edente, on aboutit \`a une contradiction avec l'hypoth\'ese
$\i\cap\i'=\{0\}$, ce qui ach\`eve de prouver le Th\'eor\`eme .\qed
\ste  Soit
$V$ un sous-espace
$\j_{0}$ invariant de
$\g$. On suppose  qu'il  est la somme de sous espaces poids de $\g$
pour
$\j_{0}$, ce qui s'\'ecrit aussi :
$$V=\sum_{\{ \lambda \in \j _{0}^{*}\vert
V^{\lambda}\not=\{0\}\}}\g^{\lambda}$$ Alors $V$ admet un
unique suppl\'ementaire $\j_{0}$-invariant, $V^{\perp}$, qui est \'egal
 \`a
la somme des sous-espaces poids de $\g$ qui ont une intersection
nulle avec $V$, soit encore :
$$V^{\perp} =\sum_{\{ \lambda \in \j _{0}^{*}\vert \g^{\lambda
}\cap V=\{0\}\}}\g^{\lambda}$$
On note $p_{V}$ (resp. $p^{V}$, la projection de $\g$ sur $V$
(resp. $V^{\perp}$) parall\`element \`a $V^{\perp}$ (resp.  $V$).
Tout sous-espace $\j_{0}$-invariant est stable sous $p^{V}$ et
$p_{V}$. \ste Si de plus $V$ est $\l$-invariant, $V^{\perp}$ est
aussi $\l$-invariant. En effet, comme $\l$ est r\'eductive dans $\g$,
$V$ admet un suppl\'ementaire $\l$-invariant qui n'est autre que
$V^{\perp}$. On voit aussi que dans ce cas, $V^{\perp}$ ne d\'epend
pas du choix de la sous-alg\`ebre de Cartan $\j_{0}$ de $\g$,
contenue dans $\l$.
On a le m\^{e}me fait pour
$\l'$ et
$\l\cap
\l'$.

\begin{theo}
Soit $B$ une forme de Manin  sur $\g$ et $\i$, $\i'$ des
sous-alg\`ebres de Lie  Lagrangiennes de $\g$,  avec
$\i$ sous $\p$ et $\i'$ sous $\p'$. On a,  gr\^{a}ce au Th\'eor\`eme 1,
$\i=\h\oplus \i_{\a}\oplus \n$, o\`u $\h=\i\cap\m$, $\i_{\a}=\i
\cap\a$. On note ${\tilde \h}=\i\cap\l$. On fait de m\^{e}me pour $\i'$.
\ste Les conditions (i) et (ii) suivantes sont \'equivalentes :
\ste (i) $(B,\i,\i')$ est un triple de Manin
\ste (ii) Notant $\i_{1}=p^{\n'}({\tilde \h }\cap\p')$,
$\i'_{1}=p^{\n}({\tilde \h' }\cap\p)$, on a :\ste
a) $\i_{1}$ et $\i'_{1}$ sont contenues dans $\l\cap\l'$, et
$(B_{1}, \i_{1}, \i'_{1}) $  est un triple de Manin dans $\l\cap
\l'$, o\`u
$B_{1}$ d\'esigne la restriction de $B$ \`a $\l\cap \l'$.\ste
b) $\n\cap \h'$ et $\n'\cap\h$ sont r\'eduits \`a z\'ero.\ste
Si l'une de ces conditions est v\'erifi\'ee, on appellera $(B_{1},
\i_{1}, \i'_{1}) $ l'ant\'ec\'edent du triple de Manin  $(B,\i,\i')$.
\end{theo}
\dem
Montrons que  (i) implique (ii). Supposons que $(B,\i,\i')$ soit un
triple de Manin dans $\g$. Pour des raisons de dimension, ceci
\'equivaut \`a
$\i\cap\i'=\{ 0\}$. Ceci implique imm\'ediatement la propri\'et\'e b) de
(ii).\ste
Montrons ensuite que $\i_{1}$ est une sous-alg\`ebre de Lie r\'eelle de
$\l\cap\l'$, et isotrope pour $B_{1}$. \ste
Etudiant les sous-espaces poids sous $\j_{0}$, on voit que :
\begin{equation}\l\cap \p' =(\l\cap\l')\oplus (\l\cap\n')
\end{equation}
Comme $\l\cap\l'$ est $\j_{0}$-invariant et que $p^{\n'}(\l\cap\n')$
est r\'eduit \`a z\'ero, on a : $$p^{\n'}(\l\cap\p')\subset \l\cap\l'$$
Il en r\'esulte que $\i_{1}$ est bien contenu dans $\l\cap\l'$. Par
ailleurs, la restriction de $p^{\n'}$ \`a $\p'$ est la projection sur
$\l'$, parall\`element \`a $\n'$. C'est donc un morphisme d'alg\`ebres de
Lie, ce qui implique que $\i_{1}$ est une sous-alg\`ebre de Lie
r\'elle de $\l\cap\l'$.\ste
Soit $X$, $X_{1}\in \i_{1}$. Ce sont des \'el\'ements de $\l\cap\l'$, et
il existe $N'$ et $N'_{1}\in \n'$ tels que $Y$ et $Y_{1}$ soient
\'el\'ements de $\i$, o\`u :
$$Y:=X+N',\>\> Y_{1}:=X_{1}+N_{1}'$$
Par ailleurs $\n'$ et $\p'$ sont orthogonaux pour $B$ (cf. la fin
de la d\'emonstration du Th\'eor\`eme 1). Un calcul imm\'ediat montre alors
que $B(Y,Y_{1})$ est \'egal \`a $B_{1}(X, X_{1})$. Comme $Y, Y_{1}\in
\i$, $B(Y, Y_{1})$ est nul. Finalement, $\i_{1}$ est isotrope pour
$B_{1}$. On montre de m\^{e}me des propri\'et\'es similaires pour
$\i_{1}^{'}$. \ste
Montrons $\i_{1}+\i_{1}' =\l\cap\l'$.
Soit $X\in \l\cap\l'$. Alors $X=I+I'$, avec $I\in \i$, $I'\in \i'$.
Ecrivons $I=H+N,\>\>I'=H'+N'$ o\`u $ H\in {\tilde \h},\>\> H'\in
{\tilde \h',\>\> N\in \n,\>\> N'\in \n'} $. On a donc :
\begin{equation}
X=H+N+H'+N'
\end{equation}
ce qui implique  :
$H=X-H'-N'-N$. On voit ainsi que $H$ est \'el\'ement de
$(\p'+\n)\cap\l$. D\'ecomposant sous l'action de $\j_{0}$, on voit
que : \begin{equation}(\p'+\n)\cap\l=\p'\cap\l\end{equation}
Finalement
$H$ est \'el\'ement de
${\tilde \h}\cap \p'$. de m\^{e}me, on voit que $H'$ est \'el\'ement de
${\tilde \h}'\cap \p$. Par  ailleurs, $\n$ et $\n'$ sont des
sous-espaces $\j_{0}$-invariants et en somme directe avec
$\l\cap\l'$. Donc, appliquant $p_{\l\cap\l'}$ \`a (2.32), on a :
$$X=p_{\l\cap\l'}(H)+p_{\l\cap\l'}(H')$$
De (2.31), on d\'eduit que la restriction de $p_{\l\cap\l'}$ \`a
$\l\cap\p'$ est \'egale \`a la restriction de $p^{\n'}$ \`a $\l\cap\p'$.
Donc, on a : $$p_{\l\cap\l'}(H)= p^{\n'}(H)\in \i_{1}$$
on obtient de m\^{e}me :
$$p_{\l\cap\l'}(H)\in \i'_{1}$$ et l'on conclut que :
$$ X\in
\i_{1}+\i'_{1}
$$
Ceci ach\`eve de prouver que : \beq\l\cap\l'=\i_{1}+\i'_{1}\eeq
Par ailleurs :$$ \l\cap\l'\>\> est
\>\>le\>\> centralisateur\>\> d\>'\>un\>\> \acute{e} l
\acute{e}ment\>\> semi-simple
\>\>de\>\>
\g,\>\> dont$$
$$\>\> l\>'\>image \>par \>\>la\>\> repr \acute{e} sentation\>\>
adjointe\>\> n\>'\>a \>\>que\>\> des\>\> valeurs
\>\>propres\>\>$$\beq
r\acute{e}
 elles\>\>\>\>\>\>\>\>\>\>\>\>\>\>\>\>\>\>\>\>\>\>\>\>\>\>\>
\>\>\>\>\>\>\>\>\>\>\>\>\>\>\>\>\>\>
\>\>\>\>\>\>\>\>\>\>\>\>\>\>\>\>
\>\>\>\>\>\>\>\>\>\>\>\>\>\>\>\>\>\>\>\>
\>\>\>\>\>\>\>\>\>\>\>\>\>\>\>\>\>\>\>\eeq
En effet
$(\l
\cap\l')\oplus ((\n'\cap\l)\oplus \n)$ est une d\'ecomposition de
Langlands d'une sous-alg\`ebre parabolique de $\g$.\ste Alors  la
restriction
$B_{1}$ de
$B$ \`a
$\l\cap\l'$ est une forme de Manin  (cf. Corollaire du Lemme 2
), et
$\i_{1}$,
$\i'_{1}$, qui sont isotropes pour $B_{1}$, sont de dimension
r\'eelles inf\'erieures ou \'egales \`a la dimension complexe de
$\l\cap\l'$.
La somme dans (2.34) est n\'ecesssairement directe, ce qui ach\`eve de
prouver que (i) implique (ii). \ste
Montrons que (ii) implique (i). Supposons satisfaites  les
conditions a) et b) de (i). Montrons que $\i\cap\i'$ est r\'eduit \`a
z\'ero. Soit $X$ un \'el\'ement de $\i\cap\i'$. Alors :
\begin{equation}X=H+N=H'+N', \>\>o\grave{u} \>\> H\in {\tilde
\h},\>\>  H'\in
{\tilde \h',\>\> N\in \n,\>\> N'\in \n'}
\end{equation}

On a alors :
$$H=H'+N'-N\in\l\cap(\p'+\n)$$ (2.33) implique que $H\in \l\cap\p'$.
De m\^{e}me, on montre que $H'\in \l'\cap\p$. Appliquant $p_{\l\cap\l'}$
\`a (2.36), on voit que :
$$p_{\l\cap\l'}(X)=p_{\l\cap\l'}(H)=p_{\l\cap\l'}(H')$$
et, gr\^{a}ce \`a la premi\`ere partie de la d\'emonstration, cela conduit
 \`a :
$$p_{\l\cap\l'}(X)=p^{\n'}(H)=p^{\n}(H')\in \i_{1}\cap\i'_{1}$$
Donc on a : $$p^{\n'}(H)=p^{\n}(H')= 0$$
Mais $p^{\n'}$ est injective sur ${\tilde \h}\cap \p'$, car
$\n'\cap{\tilde\h}$ est r\'eduit \`a z\'ero. En effet
$\n'\cap{\tilde\h}$ est contenu dans $\n'\cap\l$. On voit que cette
derni\`ere intersection est \'egal \`a $\n'\cap\m$. Donc
$\n'\cap{\tilde\h}$ est \'egal \`a $\n'\cap \h$, qui est r\'eduit \`a z\'ero,
d'apr\`es b). Donc $H$ est nul et il en va de m\^{e}me de $H'$. Alors $X$
est un \'el\'ement de $\n\cap\n'$, qui est r\'eduit \`a z\'ero, d'apr\`es nos
hypoth\`eses sur $\p$, $\p'$. Donc $X$ est nul et $\i\cap\i'$ est
r\'eduit \`a z\'ero. Alors  la somme $\i+\i'$ est directe, et l'on a
$\g=\i\oplus \i'$ pour des raisons de dimension. Ceci ach\`eve de
prouver le Th\'eor\`eme.\qed

\begin{prop}
Si $(B,\i,\i')$ est un triple de Manin sous $(\p,\p')$, d'ant\'ec\'edent
$(B_{1},\i_{1},\i'_{1})$, et si $g=nn'x\in P\cap P'$, o\`u $x\in L
\cap
L' $, $n\in N_{L'}$, $n'\in N'_{L}$, l'ant\'ec\'edent de
$(B,Ad\>g(\i),Ad\>g(\i'))$ est \'egal \`a $(B_{1},Ad\>
x(\i_{1}),Ad\>
x(\i'_{1}))$.
\end{prop}
\dem
Ecrivons ${\underline \i}=Ad\>g(\i)$ et ${\underline {\tilde
\h}}={\underline \i}\cap\l$, etc.. On note
$(B_{1},{\underline \i}_{1},{\underline \i}'_{1})$, l'ant\'ec\'edent de
$(B,Ad\>g(\i),Ad\>g(\i'))$ . On a, gr\^{a}ce \`a la Proposition 2, :
$$g=n'nx=n'x (x^{-1}n x)$$
Donc : $$Ad\>g(\i)= Ad\> n'x(\i)$$
puisque $x^{-1}n x\in N \subset I$. Comme $n'x\in (L\cap
L')N'_{L}\subset L$, cela implique : $${\underline {\tilde
\h}}= Ad\>n'x ({\tilde
\h}),$$ o\`u ${\tilde
\h}=\i\cap \l$. Mais $n'x$ est aussi \'el\'ement de $P'$. Alors, on a :
$${\underline {\tilde
\h}}\cap \p'= Ad\>n'x ({\tilde
\h}\cap \p')$$
D'o\`u l'on d\'eduit :
$${\underline \i}_{1}=p^{\n'}(Ad\>n'x ({\tilde
\h}\cap \p'))$$
Mais il est clair que la restriction de $p^{\n'}$ \`a $\p'$, n'est
autre que la projection sur $\l'$, parall\'element \`a $\n'$. Cette
restriction entrelace l'action adjointe de $P'$ sur $\p'$ avec
l'action naturelle de $P'$ sur $\l'$, identifi\'e au quotient de
$\p'$ par $\n'$ ($N'$ agit trivialement). Il en r\'esulte :
$${\underline \i}_{1}=Ad\> x(p^{\n'} ({\tilde
\h}\cap \p'))=Ad\> x(\i_{1})$$
comme d\'esir\'e. On traite de mani\`ere similaire  ${\underline
\i'}_{1}$.
\qed
\section{ Triples  de Manin pour une alg\`ebre de Lie r\'eductive
complexe : Rel\`evement}
\setcounter{equation}{0}
\begin{theo}
(i) Tout triple de Manin sous $(\p,\p')$ est conjugu\'e, par un
\'el\'ement de
$P\cap P'$ \`a un triple de Manin , $(B,\i,\i')$, sous $(\p,\p')$,
d'ant\'ec\'edent  $(B_{1},\i_{1},\i'_{1})$, satisfaisant les propi\'et\'es
suivantes :\ste
On note $\sigma$ (resp. $\sigma'$), l'af-involution  de
$\m$ (resp. $\m'$) ayant $\h$ (resp. $\h'$) pour espace de points
fixes. Il existe une \sac fondamentale ${\tilde \f}$ (resp. ${\tilde
\f'}$), de
$\i$ (resp. $\i'$), contenue dans $\i_{1}$ (resp. $\i'_{1}$) telle
que la paire $({\tilde \f}, \sigma)$ v\'erifie les conditions 1) \`a 6)
qui suivent (la paire $({\tilde \f'}, \sigma')$ v\'erifiant des
conditions similaires, num\'erot\'ees 1') \`a 6'), avec les changements
\'evidents) :\ste
1) L'application ${\sigma}$ est une af-involution  de l'alg\`ebre de
Lie $\m$, avec ensemble de points fixes $\h$ et ${\tilde \f}$ est
une \sac fondamentale de
$\i_{1}$ telle que ${\tilde \f}=\f\oplus ({\tilde \f}\cap\a)$, o\`u
$\f={\tilde \f}\cap \h$. De plus ${\tilde \f}\cap\a$ est isotrope
pour $B$, de dimension r\'eelle \'egale \`a la dimension complexe de $\a$,
et $\h$ est isotrope pour $B$. On note $\j$ le centralisateur dans
$\g$ de ${\tilde \f}$, qui est une \sac de $\g$, contenue dans
$\l\cap \l'$, d'apr\`es le Lemme 11.
\ste 2) L'intersection de $\h$ avec $\n'$ est r\'eduite \`a z\'ero.\ste
3) Il existe une sous-alg\`ebre de Borel, $\b$, de $\m$, contenant
$\j\cap\m$, et contenue dans $\m\cap\p'$, tel que
$\sigma(\b)+\b=\m$.\ste
4) Il existe une unique d\'ecomposition de Langlands de $\p_{1}$, o\`u
$\i_{1}$ est sous $\p_{1}$, $\p_{1}={\underline \l}_{1}\oplus
\n_{1}$, telle que ${\underline \l}_{1}$ contienne ${\tilde \f}$. De
plus ${\underline \m}_{1} = {\underline \l}_{1}^{der}$ est \'egal \`a
l'id\'eal d\'eriv\'e de $(\m\cap\l')\cap \sigma (\m\cap\l')$. \ste
5) Si $\alpha\in \Delta (\m, \j)$, on peut d\'efinir
${\underline \alpha}\in
\Delta (\m, \j)$ par la  condition : $\sigma(\m^{\alpha})=
\m^{{\underline \alpha}}$. On a : $\n_{1}=\oplus_{\alpha\in \Delta
(\m\cap\n', \j)\cap {\underline {\Delta (\m\cap\l',
\j)}}}\m^{{\underline \alpha}}$.
\ste 6)) La restriction de $\sigma$ \`a ${\underline \m}_{1}$ est
\'egale \`a l'af-involution dont l'espace des points fixes est \'egal
 \`a
${\underline
\h}_{1}:= \i_{1}\cap {\underline \m}_{1}$. \ste
On dit alors que le triple de Manin $(B,\i,\i')$ est li\'e \`a
$(B_{1},\i_{1},\i'_{1})$, avec lien $({\tilde \f}, {\tilde
\f'})$.\ste
(ii) R\'eciproquement si $B$ est une forme de Manin  sur $\g$, si
$B_{1}$ est sa restriction \`a $\l\cap\l'$, si
$(B_{1},\i_{1},\i'_{1})$ est un triple de Manin  pour
$\l\cap\l'$  et si
$(\sigma,{\tilde \f})$ (resp. $(\sigma',{\tilde \f'})$) v\'erifient
les propri\'et\'es 1) \`a 6) (resp. 1') \`a 6')), posant $\i={\tilde
\h}\oplus \n$, o\`u ${\tilde \h}=\h\oplus ({\tilde \f}\cap \a)$,
(resp.
$\i'={\tilde
\h'}\oplus \n'$, o\`u ${\tilde \h'}=\h'\oplus ({\tilde \f'}\cap
\a')$), alors $(B,\i,\i')$ est un triple de Manin sous
$(\p,\p')$, li\'e \`a
$(B_{1},\i_{1},\i'_{1})$, avec lien $({\tilde \f}, {\tilde
\f'})$.
\end{theo}
\begin{rem} Dans le Th\'eor\`eme, la condition 3) implique la
condition 2). En effet si 3) est satisfait, pour des raisons de
dimension,
$\b\cap\sigma (\b)$, qui contient une \sac de $\m$ ( voir (2.1)), est
r\'eduit \`a
$\j\cap\m$. Par ailleurs, comme $\b$ est contenu dans $\p'\cap\l$,
$\n'\cap\l$ est contenu dans le radical nilpotent de $\b$. Par
suite, on a :
$$(\n'\cap\l)\cap \sigma (\n'\cap\l)\subset (\n'\cap\l)\cap
\j =\{0\}$$
Donc, $\n'\cap\l\cap \h$ est r\'eduit \`a z\'ero. Comme $\h$ est contenu
dans $\l$, $\n'\cap\h$, $\n'\cap\l$ est \'egalement r\'eduit \`a z\'ero,
comme d\'esir\'e.
\end{rem}
\dem
D\'emontrons (i).
Soit $(B,{\underline \i},{\underline \i}')$ un triple de Manin pour
$\g$, sous
$(\p,\p')$. On note ${\underline \h}={\underline \i}\cap \m$,
etc. Comme ${\underline \i}+{\underline \i'}= \g$, on a :
$${\underline \i}+\p'=\g$$ Appliquant $p_{\l}$ \`a cette \'egalit\'e, on
en d\'eduit  :
$$({\underline \i}\cap \l)+(\p'\cap\l)=\l$$
On applique encore la projection de $\l$ sur $\m$, parall\`element \`a
$\a$ pour obtenir :$${\underline \h}+(\p'\cap \m)=\m$$
En cons\'equence, ${\underline H }(P'\cap M)^{0}$ est ouvert dans $M$.
Or
$(P'\cap M)^{0}$ est le sous groupe parabolique de $M$, d'alg\`ebre
de Lie
$\p'\cap\m$. Par ailleurs $\sigma$ \'etant une af-involution, $\m$
est le produit d'id\'eaux $\m_{j}$, invariants par $\sigma$ et sur
lesquels induit :
\ste soit une conjugaison par rapport \`a une forme r\'eelle,\ste  soit
''l' \'echange des facteurs ''  de deux id\'eaux
isomorphes dont $\m_{j}$
est la somme.\ste Il r\'esulte alors de [M2], [M1], que
${\underline \h}\cap\p'$ contient une \sac fondamentale ${\underline
\f}$ de
${\underline \h}$ et une sous-alg\`ebre de Borel de $\m$, contenant
${\underline
\f}$, contenue dans $\p'\cap\m$, et telle que  :
\begin{equation}
\sigma({\underline \b})+{\underline \b}=\m
\end{equation}
D'apr\`es le Lemme 11, le centralisateur  ${\underline
\j}$ de ${\tilde {\underline \f}}:={\underline \f}+{\underline
\i}_{\a}$ dans $\g$ est une \sac  de $\g$, contenue dans
$\l$. De la d\'efinition des af-involutions, il r\'esulte que toute \sac
de
${\underline \h}$ contient des \'el\'ements r\'eguliers de $\m$. Il
r\'esulte alors de [Bor], Proposition 11.15, que ${\underline
\j}\cap\m$ est contenu dans $\b$. Donc
${\underline \j}=({\underline \j}\cap\m)\oplus
\a$ est contenu dans $\p'\cap\l$. C'est une \sac de $\p'\cap\l$,
donc elle est conjugu\'ee \`a $\j_{0}$, par un \'el\'ement du sous-groupe
analytique de $G$, d'alg\`ebre de Lie $\p'\cap\l$. Mais, d'apr\`es
(2.31), on a  : $\p'\cap\l=(\l\cap\l')\oplus (\n'\cap\l)$ et ce
sous-groupe analytique est \'egal \`a $(L\cap L')N'_{L}$, puisque
$L\cap L'$ est connexe, d'apr\`es la Proposition 2.\ste
Donc, il existe $n'\in N'_{L}$, $x\in L\cap L'$, tels que
$$Ad\>xn'({\underline
\j})=\j_{0}$$
soit encore :
\begin{equation} Ad\> n'({\underline \j})=Ad\> x^{-1}(\j_{0})\subset
\l\cap
\l'
\end{equation}
On trouve de m\^{e}me ${\tilde {\underline \f'}}$, $ {\underline
\b'}$, $ {\underline \j'}$ et $x'\in L\cap L'$, $n\in N_{L'}$,
v\'erifiant des propri\'et\'es similaires. On pose :
$$u=nn', \>\> \i=Ad\> u({\underline \i}), \>\> \i'=Ad\> u({\underline
\i'}) $$
Comme $n$ et $n'$ commutent et que $\n$ est un id\'eal de $\i$, on a :
$$Ad\> u({\underline \i})=Ad\> n'({\underline \i})$$
et de m\^{e}me :
$$Ad\> u({\underline \i'})=Ad\> n({\underline \i'})$$
On pose alors :\beq {\tilde \f}= Ad\> n'({\tilde {\underline \f}}),
\>\> {\tilde \f'}= Ad\> n({\tilde {\underline \f}}'), \>\> \b= Ad\>
n'({\underline \b}) , \>\> \b'= Ad\>
n({\underline \b'})
\eeq
On voit alors que $(B,\i,\i')$ est un triple de Manin, conjugu\'e par
$u$ \`a
$(B,{\underline \i},{\underline \i}')$ et   sous
$(\p,\p')$. \ste
On va voir que ${\tilde \f}$ a les propri\'et\'es voulues.
D'abord, ${\tilde {\underline \f}}$ est une \sac fondamentale de
${\underline \i}$, d'apr\`es le Lemme 11. Par conjugaison, on en
d\'eduit que ${\tilde \f}$ est une \sac fondamentale de $\i$. Le
centralisateur,  $\j$, de ${\tilde \f}$ v\'erifie  \beq
\j= Ad\> n'({\underline \j})
\eeq
donc est contenu dans $\l\cap\l'$, d'apr\`es (3.2). Alors, d'apr\`es
le Lemme 11, on a bien
${\tilde \f}=\f\oplus ({\tilde \f}\cap\a)$, o\`u
$\f={\tilde \f}\cap \h$, et ${\tilde \f}\cap\a$ est \'egal \`a
$\i_{\a}$, donc isotrope pour $B$ et de la dimension voulue,
d'apr\`es
le Th\'eor\`eme 1. De m\^{e}me $\h$ est isotrope pour $B$.\ste On a vu que
${\tilde \f}$ est contenu dans $\i\cap\l\cap\l'$, donc dans
${\tilde
\h}\cap\p'$. De plus
$p^{\n'}$ est l'identit\'e sur $\l'$. Donc $ {\tilde \f}$ est contenu
dans $\i_{1}$. Par ailleurs, comme ${\tilde \f}$ est une \sac
de $\i$, contenue dans ${\tilde \h}\cap\p'$, c'est  une \sac de
${\tilde \h}\cap\p'$ (cf. [Bou], Ch. VII, Paragraphe 2.1, Exemple
3), et par projection , c'est une \sac de $\i_{1}$ (cf. l.c.,
Corollaire 2 de la Proposition 4). \ste Il reste \`a voir, pour
achever de v\'erifier 1), que cette \sac de $\i_{1}$ est
fondamentale.\ste
 D'apr\`es le Lemme 11 (iv), il existe une unique d\'ecomposition de
Langlands $\p_{1}={\underline \l_{1}}\oplus\n_{1}$, telle que
${\underline \l_{1}}$ contienne ${\tilde \f}$. On note ${\underline
\m_{1}}={\underline \l_{1}}^{der}$. Il suffit de voir que
${\tilde
\f}\cap {\underline
\m_{1}}$ est une \sac fondamentale de ${\underline
\h_{1}}$. Pour cela, il suffit de voir qu'aucune racine de $
{\tilde \f}\cap {\underline
\m_{1}}$ dans ${\underline
\m_{1}}$ n'est r\'eelle. D'apr\`es le Lemme 11 (iv), ${\tilde \f}=
({\tilde \f}\cap {\underline
\m_{1}})\oplus ({\tilde \f}\cap {\underline
\a_{1}})$, o\`u ${\underline
\a_{1}}$ est le centre de ${\underline
\l_{1}}$. Alors, une racine $\alpha $ de $
{\tilde \f}\cap {\underline
\m_{1}}$ dans ${\underline
\m_{1}}$, prolong\'ee par z\'ero sur ${\tilde \f}\cap {\underline
\a_{1}}$ est une racine de ${\tilde \f}$ dans $\m$. Mais alors,
comme $\f$ est une \sac fondamentale de $\h$, $\alpha $ n'est pas
r\'eelle sur $\f$. Ceci  prouve que ${\tilde \f}$ est une \sac
fondamentale de
$\i_{1}$.\ste
2) r\'esulte du fait que l'intersection de $\i$ et $\i'$
est r\'eduite \`a z\'ero, car $\i'$ contient $\n'$.
\ste L' \'equation (3.1) assure, par transport de structure que
3) est satisfait.
Comme ${\tilde \f}$ est une \sac de $\i_{1}$, l'existence de
${\underline
\l_{1}}$ contenant
${\tilde \f}$, r\'esulte du Lemme 11 (iv).
\ste D\'ecrivons plus pr\'ecis\'ement $\i_{1}$. Pour cela commencons par
d\'ecrire l'action de $\f$ sur ${\tilde \h}\cap \p'$.
Soit $\alpha\in \Delta (\m,\j)$.  Comme les \'el\'ements de $\f$ sont
fix\'es  par
$\sigma$, $\j\cap\m$, qui est le centralisateur de $\f$ dans $\m$,
est $\sigma$-invariant. Il r\'esulte du Corollaire du Lemme 6, que
$\m^{\alpha}$ est alors contenu dans un id\'eal de $\m$, stable sur
lequel $\sigma$ est soit lin\'eaire soit antilin\'eaire
Donc,
$\sigma$ est soit lin\'eaire soit antilin\'eaire sur $\m^{\alpha}$.
Comme, pour $X\in
\m^{\alpha}$ et $H\in \j $, on a :
\beq [H,\sigma(X)]=[\sigma\sigma (H), \sigma (X)]= \sigma
([\sigma (H),X])=\sigma(\alpha(\sigma(H))X),\eeq et on peut
d\'efinir
${\underline \alpha}\in  \Delta(\m,\j)$ par :
l'\'egalit\'e
\begin{equation}
\sigma (\m^{\alpha})= \m^{{\underline \alpha}}
\end{equation}
Par ailleurs $\b$ et $\sigma (\b)$ sont des
 sous-alg\`ebres de Borel
oppos\'ees relativement \`a $\j\cap\m$, car 3) a \'et\'e v\'erifi\'e. On
d\'eduit
de ce qui pr\'ec\`ede :
\begin{equation}
\alpha\not= {\underline \alpha},\>\> \alpha\in
\Delta (\m,\f)
\end{equation}
On note :
\begin{equation}
\h_{\alpha}:= (\m^{\alpha}+\m^{{\underline \alpha}})^{\sigma},\>\>
\alpha\in \Delta (\m,\f)
\end{equation}
On a :
\begin{equation}
\h=\f\oplus _{\alpha\in \Delta (\b,\j)}\h_{\alpha}
\end{equation}
On v\'erifie ais\'ement  que $\h_{\alpha}$ est invariant sous $\f$.
Montrons :
\beq
\h_{\alpha}\>\> est\>\> irr\acute{e} ductible\>\> sous \>\>\f
\eeq
Comme $\alpha$ est diff\'erent de ${\underline \alpha}$,
$\h_{\alpha}$, qui est \'egal \`a $\{X+\sigma(X)\vert X\in
\m^{\alpha}\}$, est de dimension deux sur $\R$. S'il a un
sous-module irr\'eductible non nul, $\R X$, $X$ doit se
transformer sous $\f$ par un caract\`ere r\'eel de $\f$. Celui-ci doit
\^{e}tre  \'egal \`a la restriction de $\alpha$ ou ${\underline
\alpha}$ \`a $\f$. Alors (3.10) va r\'esulter de la d\'emonstration de :
\beq
Pour \>\>tout\>\> \alpha \in \Delta(\m,\j),\>\> \alpha
\>\>restreint \>\>\grave{a} \>\>\f \>\>n\>'\>est\>\> pas\>\>
r\acute{e} elle
\eeq
Distinguons deux cas.\ste
a) $\m^{\alpha}$ est contenu dans un id\'eal $\m_{0}$, stable par
$\sigma$, sur lequel $\sigma$ est antilin\'eaire. \ste Pour d\'emontrer
(3.11), on se r\'eduit au cas o\`u $\m=\m_{0}$, i.e. on peut supposer
$\sigma$ antilin\'eaire. (3.5) montre que la
restriction de
${\underline \alpha}$ \`a $\f$ est \'egale \`a la conjugaison complexe
de la restriction de $\alpha$ \`a $\f$, donc $\alpha$ est r\'eelle
sur $\f$ si et seulement si ces deux restrictions sont \'egales.
Mais $\h$ \'etant, dans notre cas , une forme r\'eelle de $\m$, $\f$
est une forme r\'eelle de $\j$.  Si les deux restrictions \`a $\f$
\'etaient \'egales , on aurait alors $\alpha={\underline \alpha}$,
ce qui
n'est pas. (3.11) en r\'esulte dans ce cas.\ste
b) $\m^{\alpha}$ est contenu dans un id\'eal de $\m$, produit de
deux id\'eaux simples, $\m'$, $\m''$, et $\sigma$ restreint \`a
$\m'\times \m''$ est $\C$ lin\'eaire et de la forme :
$$(X', X'')\mapsto (\tau^{-1}(X''), \tau(X') ),\>\> (X',
X'')\in \m'\times \m''$$
o\`u $\tau $ est un isomorphisme $\C$-lin\'eaire entre les alg\`ebres
de Lie $\m'$ et $\m''$. Ici aussi, on se r\'eduit \`a prouver
l'assertion lorsque $\m=\m'\times\m''$. L '\'etude des points
fixes de
$\sigma$, montre ais\'ement que, dans ce cas, $\f$ est un
sous-espace vectoriel complexe de $\j$, sur lequel la
restriction de $\alpha $ est non nulle , et automatiquement
$\C$-lin\'eaire, donc n'est pas r\'eelle. Ce qui ach\`eve de  prouver
(3.11).\ste
En d\'ecomposant l'action de $\j_{0}$ sur $\l\cap \p'$, on trouve :
$\l\cap\p'=(\m\cap\p')\oplus \a$, puis, tenant compte de la
d\'efinition de ${\tilde \h}$ :
${\tilde \h}\cap \p'= (\h\cap\p') \oplus \i_{\a}$
D\'ecomposant $ \h\cap \p'$ en repr\'esentations irr\'eductibles
sous $\f$, on voit que : $$ \h\cap \p'=\f \oplus
 _{\alpha\in
\Delta(\b,\j),\>\> \h_{\alpha}\subset\p'}\h_{\alpha}$$
Il est clair que $\h_{\alpha}$ est contenu dans $\p'$ si et
seulement si c'est le cas de $\m^{\alpha}$ et de $\m^{\underline
\alpha}$. Mais on a :
$$\m\cap\p'=(\m\cap\l')\oplus (\m\cap\n')
$$
De plus comme $\b$ est contenu dans $\p'$, $\m\cap\n'$
est contenu dans $\b$, c'est \`a \ste dire :
\beq \m\cap\n'=\b\cap \n'
\eeq En particulier, si $\alpha \in \Delta
(\b,\j)$, ${\underline \alpha}$ n'est pas \'el\'ement de $\Delta
(\m\cap\n', \j)$. Par suite :

$$
\h\cap\p'=\f \oplus (\oplus_{\alpha\in \Delta (\b\cap \l',\j)
\cap
  {\underline {\Delta (\m\cap \l',\j)}}}\h_{\alpha})\oplus(
\oplus_{\alpha\in \Delta (\b\cap \n',\j)
\cap
  {\underline {\Delta (\m\cap \l',\j)}}}\h_{\alpha})
$$

Etudions l'image, $\i_{1}$, de ${\tilde \h}\cap \p'$ par $p^{\n'}$.
 Si
$\alpha\in \Delta (\b\cap \l',\j)
\cap
  {\underline {\Delta (\m\cap \l',\j)}}$, $\h_{\alpha}$ est contenu
dans $\m\cap\l'$, et  $p^{\n'}(\h_{\alpha})$ est \'egal \`a
$\h_{\alpha}$. Si $\alpha\in \Delta (\b\cap \n',\j )
\cap
  {\underline {\Delta (\m\cap \l',\j)}}$, $\m^{\alpha}$ est contenu
dans $\n'$ et $\m^{{\underline \alpha}}$ est contenu dans $\l'$.
Comme $\h_{\alpha}=\{X+\sigma(X)\vert X\in \m^{\alpha}\}$,
$p^{\n'}(\h_{\alpha})$ est \'egal \`a
$\m^{{\underline \alpha}}$. Enfin $p^{\n'}$ est l'identit\'e sur
${\tilde \f}=\f\oplus \i_{\a}$.\ste   Finalement :
\beq
p^{\n'}({\tilde \h}\cap \p')=\u_{1}\oplus \v_{1}\oplus \i_{\a}
\eeq
o\`u :
\beq
\u_{1}=\f\oplus(\oplus_{\alpha\in \Delta (\b\cap \l',\j)
\cap
  {\underline {\Delta (\m\cap \l',\j)}}}\h_{\alpha})
\eeq
et \beq
\v_{1}=\oplus_{\alpha\in \Delta (\b\cap \n',\j)
\cap
  {\underline {\Delta (\m\cap \l',\j)}}}\m^{{\underline \alpha}}
\eeq
On remarque que :
\beq
\u_{1}= ((\m\cap \l')\cap \sigma (\m\cap\l'))^{\sigma}
\eeq
Donc $\u_{1} $ est une sous-alg\`ebre de Lie de $\m\cap\l'$,
r\'eductive. On note $\z_{1}$ son centre. Alors : $$\r_{1}:=
\z_{1}\oplus
\v_{1}\oplus \i_{\a}$$ est un id\'eal r\'esoluble de $\i_{1}=
p^{\n'}({\tilde \h}\cap \p')$, admettant la sous alg\`ebre de Lie
semi-simple $\u_{1}^{der}$ comme suppl\'ementaire dans $\i_{1}$. Donc
$\r_{1}$ est le radical de $\u_{1}$. Etudiant l'id\'eal d\'eriv\'e,
$\i_{1}^{der}$ de
$\i_{1}$, on a :
$$\i_{1}^{der}\subset \u_{1}^{der}\oplus \v_{1}$$
Par ailleurs $\u_{1}^{der}$ est contenu dans $\i_{1}^{der}$, de
m\^{e}me que $\v_{1}$, car $[\f,\v_{1}]=\v_{1}$, d'apr\`es (3.11).
Finalement, on obtient :
$$\i_{1}^{der}= \u_{1}^{der}\oplus \v_{1}, \>\>
\r_{1}\cap\i_{1}^{der}=\v_{1}$$
Donc $\v_{1}$ est \'egal au radical nilpotent de $\i_{1}$ (cf.
[Bou], Ch. I, Paragraphe 6.4, Proposition 6), i.e. :
\beq
\v_{1}=\n_{1}
\eeq
La propri\'et\'e 5) est prouv\'ee.\ste
Prouvons 4).
On a vu que le centralisateur  $\j$ de ${\tilde \f}$ est une \sac de
$\l$,  \'egale \`a $(\j\cap\m)+\a$. Par
ailleurs, comme
${\tilde
\f}$ est contenue dans $\i_{1}$, donc dans $\l'$, ce centralisateur
contient
$\a'$, puis est contenu dans $\l'$. Donc, $\j$ est une \sac
de $\l\cap\l'$. On pose :
\beq
{\tilde {\underline \l_{1}}}:= ((\m\cap \l')\cap\sigma (\m\cap\l'))
\oplus \a
\eeq
qui contient $\j$.
 On pose  ${\underline  \p_{1}}:={\tilde {\underline \l_{1}}}\oplus
\v_{1}$, qui contient  la
sous-alg\`ebre de Borel de $\l\cap\l'$, $\l'\cap (\sigma(\b)\oplus
\a)$. \ste
 On  va voir que  ${\underline  \p_{1}}$ est une sous-alg\`ebre de
Lie   de
$\l\cap\l'$. L'\'equation (3.12) jointe \`a la d\'efinition de
$\u_{1}$ et
$\v_{1}$ (cf (3.14), (3.15) et (3.16)), permet de voir que le
crochet de
$\u_{1}$ et $\v_{1}$ est contenu dans $\v_{1}$. Ceci implique
que ${\underline  \p_{1}}$ est une sous alg\`ebre de Lie de
$\l\cap\l'$. C'est une sous-alg\`ebre parabolique, car elle
contient une sous-alg\`ebre de Borel.\ste L'analyse des poids sous
$\j$ montre que
${{\tilde {\underline \l_{1}}}}\oplus\v_{1}$
est une d\'ecomposition de Langlands de la sous-alg\`ebre
parabolique, ${\underline  \p_{1}}$, de $\l\cap\l'$. En particulier
${\underline  \p_{1}}$ admet $\v_{1}$ comme radical
nilpotent. Donc $\p_{1}=  {\underline  \p_{1}}$, d'apr\`es (3.17)
et (1.17), et
${\tilde {\underline \l_{1}}}\oplus
\n_{1}$ est une d\'ecomposition de Langlands de $\p_{1}$, telle que
${\tilde {\underline \l_{1}}}$ contienne ${\tilde \f}$. C'est la
seule, d'apr\`es le Lemme 11 (iv), appliqu\'e \`a $\i_{1}$. Donc
${\tilde {\underline \l_{1}}}={\underline \l_{1}}$\ste  De plus,
comme
$\a$ est central dans
$\l$, il r\'esulte  (3.18) que ${\underline  \m_{1}}$ est l'id\'eal
d\'eriv\'e de $(\m\cap \l')\cap(\sigma
(\m\cap\l')$. D'apr\`es le Th\'eor\`eme 1, on sait que ${\underline
\h_{1}}:= \i_{1}\cap {\underline  \m_{1}}$ est l'espace des points
fixes d'une af-involution de ${\underline  \m_{1}}$ que l'on note
${\underline \sigma}_{1}$.
${\underline  \m_{1}}$.\ste
 En outre, on a ${\underline
\h_{1}}=(\i_{1}\cap {\underline  \l_{1}})^{der}$, d'o\`u l'on d\'eduit
:
\beq
{\underline  \h_{1}}=\u_{1}^{der}
\eeq On sait, d'apr\`es (3.16), que $\u_{1}$ est l'espace des points
fixes de la restriction de $\sigma$ \`a ${\tilde {\underline
\m_{1}}}$, dont l'id\'eal d\'eriv\'e est \'egal \`a ${\underline
\m_{1}}$.Comme ${\tilde {\underline
\m_{1}}}$ est r\'eductive,  passant \`a
l'id\'eal d\'eriv\'e, on voit que ${\underline
\h_{1}}=\u_{1}^{der}$ est \'egal \`a l'id\'eal d\'eriv\'e de
l'ensemble ${\underline \h_{1}}'$ des points fixes de la restriction
${\underline
\sigma_{1}}'$, de
$\sigma$ \`a
${\underline
\m_{1}}$. Montrons que ${\underline \h_{1}}'={\underline
\h_{1}}$. Notons $\q_{1}$ (resp. $\q'_{1}$), l' orthogonal de
${\underline \h_{1}}$ (resp. ${\underline \h_{1}}'$) pour la forme
de Killing de l'alg\`ebre de Lie  ${\underline  \m_{1}}$, regard\'ee
comme alg\`ebre de Lie r\'eelle. Alors, on a :
$${\underline
\m_{1}}={\underline \h_{1}}\oplus\q_{1}={\underline
\h_{1}}'\oplus{\q_{1}}'$$
D'o\`u l'on d\'eduit :
$${\underline \m_{1}}^{der}={{\underline \h'_{1}}}^{der}+
[\q'_{1},\q'_{1}]+ [{\underline \h'_{1}},\q'_{1}]
$$
Mais $\q'_{1}$ est contenu dans $\q_{1}$ et $[\q_{1},\q_{1}]$ est
contenu dans ${\underline \h_{1}}$. D'autre part $[{\underline
\h'_{1}},\q'_{1}]$ est contenu dans $\q'_{1}$. Finalement, tenant
compte du fait que ${\underline
\m_{1}}$ est semi-simple, on a :
$${\underline
\m_{1}}\subset {\underline \h_{1}}+\q'_{1}$$
D'o\`u l'on d\'eduit l'\'egalit\'e de ${\underline \h_{1}}$ et
${\underline \h_{1}}'$, et par suite celle de ${\underline
\sigma}_{1}$ et ${\underline \sigma}'_{1}$. Ceci ach\`eve de
prouver 5) et 6).\ste  On vient d'achever la preuve de (i).
\ste Passons \`a la preuve de (ii), dont on retient les hypoth\`eses.
D'abord $\i$ et $\i'$ sont des sous-alg\`ebres de Lie r\'eelles,
isotropes,  de $\g$, dont la dimension r\'eelle est \'egale \`a la
dimension complexe de $\g$, d'apr\`es la condition 1) et le
Th\'eor\`eme 1. On va appliquer le Th\'eor\`eme 3 pour voir que
$(B,\i,\i')$ est un triple de Manin. \ste
Une analyse de l'action de $\f$, similaire \`a celle qui a \'et\'e
faite au cours de la preuve de (i) conduit \`a  l'analogue des
\'equations (3.13)  \`a (3.16). L'analogue de l'\'equation  (3.16) et la
condition 6) montrent  que
$\u_{1}^{der}$ contient ${\underline \h_{1}}:={\underline
\m_{1}}\cap\i_{1}$.\ste  Par ailleurs
${\tilde
\f}$ est contenu dans $p^{\n'}({\tilde \h}\cap \p')$ , car ${\tilde
\f}\subset \l\cap\l'\cap {\tilde \h}$, d'apr\`es 1).\ste D'autre part,
d'apr\`es 5) : $
\v_{1}=\n_{1}$.\ste
Finalement $p^{\n'}({\tilde \h}\cap \p')$ contient $\i_{1}$. Par
ailleurs, comme $\i$ est isotrope pour $B$, le raisonnement du d\'ebut
du Th\'eor\`eme 2, montre que $p^{\n'}({\tilde \h}\cap
\p')$ est isotrope pour $B_{1}$. Donc, pour des raisons de
dimension, on a :
$$p^{\n'}({\tilde \h}\cap \p')=\i_{1}$$
De m\^{e}me, on voit que :
$$p^{\n}({\tilde \h'}
\cap \p)=\i'_{1}$$
Tenant compte de 3), le Th\'eor\`eme 3 montre que $(B,\i,\i')$ est un
triple de Manin pour $\g$, d'ant\'ec\'edent $(B_{1},\i_{1},\i'_{1})$.
\ste Montrons que ${\tilde \f}$ est une \sac fondamentale de $\i$.
D'apr\`es (3.11), aucune racine de $\f$ dans $\m$ n'est r\'eelle, donc
$\f$ est une
\sac fondamentale de $\h$. Alors
${\tilde
\f}$ est une \sac fondamentale de $\i$, d'apr\`es le Lemme . On a
le m\^{e}me r\'esultat pour ${\tilde
\f'}$. Ceci ach\`eve de prouver le Th\'eor\`eme. \qed
\begin{prop} (i) Tout triple de Manin sous $(\p,\p')$ est conjugu\'e,
par un \'el\'ement de $G$, \`a un triple de Manin, $(B,\i,\i')$,  sous
$(\p,\p')$, pour lequel il existe :\ste
1) Une suite finie strictement d\'ecroissante de sous-alg\`ebres de
Lie r\'eductives complexes de $\g$, $\g_{0}=\g$, $\g_{1}=\l\cap
\l'$,...,$\g_{k}$,..., contenant $\j_{0}$, se terminant \`a
$\g_{p_{0}}=\j_{0}$. \ste  2) Une suite de  couples de
sous-alg\`ebres paraboliques des
$\g_{k}$,  $(\p_{k},\p'_{k})$, $\p_{k}$
contenant
$\b_{0}\cap
\g_{k}$, $\p'_{k}$ contenant $\b'_{0}\cap
\g_{k}$, avec $\p_{0}=\p$, $\p'_{0}=\p'$, et telle que, notant
$\l_{k}\oplus \n_{k}$ (resp. $\l'_{k}\oplus \n'_{k}$ ) la
d\'ecomposition de Langlands de $\p_{k}$ (resp. $\p'_{k}$) avec
$\j_{0}\subset \l_{k}, \>\l'_{k}$, on a :
$$\g_{k}=\l_{k-1}\cap \l'_{k-1},\>\>i=1,\dots , p_{0}$$
3) Une suite de triples de Manin,  $(B_{k},\i_{k},\i'_{k})$,  dans
$\g_{k}$,
$k=0,
\dots , p_{0}$, o\`u $(B_{0},\i_{0},\i'_{0})= (B,\i,\i')$, chacun
\'etant li\'e \`a son ant\'ec\'edent. On notera $({\tilde \f}_{k},{\tilde
\f'}_{k}) $ un lien entre $(B_{k},\i_{k},\i'_{k})$ et
$(B_{k+1},\i_{k+1},\i'_{k+1})$, $k=0, \dots , p_{0}-1$.\ste   On
appelera une telle donn\'ee une tour standard de triples de Manin dans
$\g$. On dira aussi que $(B,\i,\i')$ est fortement standard.\ste
(ii) La suite des $(\g_{k},\p_{k},\p'_{k}, B_{k})$ ne d\'epend que de
la classe de conjugaison sous $G$ du triple de Manin initial. En
particulier
$p_{0}$ ne d\'epend pas de la tour choisie. On l'appelle la hauteur
(de la classe de conjugaison sous $G$) du triple de Manin de
d\'epart. De plus, si deux tours standards,
$(B_{k},\i_{k},\i'_{k})$, $(B_{k},{\underline
\i_{k}},{\underline\i'_{k}})$, sont associ\'ees \`a un m\^{e}me triple de
Manin, pour tout
$k=1,\dots, p_{0}$, il existe $g_{k}\in G_{k}$  conjuguant
$\i_{k}$ et ${\underline
\i_{k}}$ et aussi  $\i'_{k}$ et ${\underline
\i'_{k}}$.\ste
(iii)  Le triple de Manin dans $\j_{0}$,
$(B_{p_{0}},\i_{p_{0}},\i'_{p_{0}})$, ne d\'epend que (de la classe de
conjugaison sous $G$) du triple de Manin de d\'epart. En particulier
il ne d\'epend pas de la tour standard associ\'ee \`a ce triple de Manin.
On l'appelle le socle de la classe de conjugaison sous $G$ de
triples de Manin consid\'er\'ee. Notons que $(\i_{p_{0}},\i'_{p_{0}})$
est juste une paire de sous-espaces vectoriels r\'eels   de $\j_{0}$,
isotropes pour la restriction, $B_{p_{0}}$, de $B$, de dimensions
maximales et en somme directe.\ste
(iv) Il existe $y_{k}\in I_{k}$, $k=1,\dots , p_{0}$ tel que :
$${\tilde \f}_{k}= Ad \>y_{k+1}\dots Ad
\>y_{{p_{0}}}(\i_{p_{0}})\subset \i_{k}\cap\i_{k+1},\>\> k=0,
\dots  ,p_{0}-1$$
On notera $x_{1}=y_{1}\dots y_{p_{0}}\in I_{1}\dots
I_{p_{0}}\subset L\cap L' $, de sorte que ${\tilde
\f}_{0}=Ad\>x_{1}(\i_{p_{0}})$. On introduit un \'el\'ement $x'_{1}$
ayant les m\^{e}mes propri\'et\'es pour ${\tilde
\f'}_{0}$ et $\i'_{p_{0}}$
\end{prop}
\dem
Montrons (i) par r\'ecurrence sur le rang de
$\g^{der}$. Le cas o\`u $\g^{der}$ est nul est clair. Suposons le
non nul.\ste D'apr\`es le Th\`eor\`eme 3 (i), \'etant donn\'e un triple de
Manin sous
$(\p,\p')$, il existe un triple de Manin, $(B,{\underline
\i},{\underline\i'})$,  qui lui est conjugu\'e sous $G$, et qui est
li\'e \ste \`a son ant\'ec\'edent $(B_{1},{\underline
\i_{1}},{\underline\i'_{1}})$. Comme  $\l\cap\l'$ est diff\'erent de
$\g$, d'apr\`es  le Th\'eor\`eme 2,
l'hypoth\'ese de r\'ecurrence implique qu' il existe
$g_{1}\in L\cap L'= ( L\cap L')^{0}$, tel que
$(B_{1},\i_{1},\i'_{1})=(B_{1},Ad\> g_{1}({\underline
\i_{1}}),Ad\> g_{1}({\underline\i'_{1}})) $ soit le premier \'el\'ement
d'une tour standard,  $T_{1}$, de $\g_{1}:=\l\cap\l'$. Mais il est
facile de voir, par transport de structure et gr\^{a}ce \`a la
Propostion 4, que
$(B,\i,\i'):=(B,Ad\> g_{1}({\underline
\i}),Ad\> g_{1}({\underline\i'}))$ est li\'e \`a
$(B_{1},\i_{1},\i'_{1})$, un lien entre les deux \'etant obtenu par
conjugaison par $g_{1}$ d'un lien entre $(B,{\underline
\i},{\underline\i'})$ et $(B_{1},{\underline
\i_{1}},{\underline\i'_{1}})$. Adjoignant $(B,\i,\i')$ \`a la tour
standard, $T_{1}$, dans $\g_{1}$, on obtient une tour standard
dans $\g$, qui poss\`ede les propri\'et\'es voulues. Ceci prouve (i).\ste
Montrons (ii) par r\'ecurrence sur le rang de
$\g^{der}$. Si deux tours standards,
$(B_{k},\i_{k},\i'_{k})$, $(B_{k},{\underline
\i_{k}},{\underline\i'_{k}})$, sont associ\'ees \`a un
m\^{e}me triple de Manin sous $(\p,\p')$, on a $\p_{0}={\underline
\p_{0}}=\p$,
$\p'_{0}={\underline \p'_{0}}=\p$, d'apr\`es la Proposition 1, donc
$\g_{1}={\underline
\g_{1}}$. De plus, d'apr\`es la Proposition 2 (i),
$(B_{0},\i_{0},\i'_{0})$,
$(B_{0},{\underline
\i_{0}},{\underline\i'_{0}})$ sont conjugu\'es par un \'el\'ement de
$P\cap P'$.   Alors, d'apr\`es la Proposition 4, les ant\'ec\'edents de
$(B_{0},
\i_{0},\i'_{0})$,
$(B_{0},{\underline
\i_{0}},{\underline\i'_{0}})$ sont conjugu\'es par un \'el\'ement de
$G_{1}=L\cap L'$. On applique l'hypoth\`ese de r\'ecurrence pour
achever de prouver (ii).\ste
(iii) r\'esulte de (ii) appliqu\'e \`a l'indice $k=p_{0}$. \ste
Pour (iv), on proc\`ede par r\'ecurrence descendante sur $k$. Il suffit
alors d'obser-\ste ver que, ${\tilde \f}_{k}$ et ${\tilde \f}_{k+1}$
\'etant deux
sous-alg\`ebres de Cartan fondamentales de $\i_{k}$, d'apr\`es la
d\'efinition des liens (cf. Th\'eor\`eme 4), elles sont  conjugu\'ees par
un \'el\'ement de
$I_{k+1}$.\qed\ste
On s'int\'eresse maintenant \`a l'utilisation de la partie
r\'eciproque du Th\'eor\`eme 4. On se fixe une forme de Manin  $B$
sur $\g$. On note
$B_{1}$ sa restriction \`a $\l\cap \l'$. On
se donne un triple de Manin, $(B_{1},\i_{1},\i'_{1})$,   fortement
standard dans $\l\cap\l'$. On emploie les
notations de la Proposition 5, pour une tour standard d\'ebutant par
ce triple de Manin. On fixe  ${\tilde \f}$ (resp. ${\tilde \f'}$)
une
\sac fondamentale de $\i_{1}$ (resp. $\i'_{1}$).  Proc\`edant comme
dans la preuve du point (iv) de la Proposition 5, on trouve
$x_{1}\in I_{1}\dots I_{p_{0}},\>\> x'_{1}\in I'_{1}\dots
I'_{p_{0}}$, avec ${\tilde
\f}=Ad\>x_{1}(\i_{p_{0}})$, ${\tilde
\f'}=Ad\>x'_{1}(\i'_{p_{0}})$.
Si $\sigma$ est une af-involution  de $\m$, on notera
${\overline \sigma}= Ad\>{x_{1}}^{-1}\circ\sigma
\circ Ad\> {x_{1}}_{\vert \m}$ et ${\overline \h}= \m^{{\overline
\sigma}}$, etc.\ste
L'int\'er\^{e}t du Lemme suivant  r\'eside dans le fait qu'il r\'eduit une
partie de la construction d'un rel\`evement de
$(B_{1},\i_{1},\i'_{1})$, avec lien $({\tilde
\f},{\tilde
\f'})$, \`a des questions sur le socle de celui ci.

\begin{lem}
Soit $\sigma$ une af-involution
l'alg\`ebre de Lie $\m$.\ste La paire $(\sigma, {\tilde \f})$ v\'erifie
les conditions 1) \`a 5) du Th\'eor\`eme 4, si et seulement si
$({\overline \sigma}, \i_{p_{0}})$  v\'erifie des
conditions similaires, o\`u l'on remplace ${\underline \l_{1}}$ par
$\l_{1}$ (o\`u $\l_{1}\oplus \n_{1}$ est la d\'ecomposition de
Langlands de $\p_{1}$ telle que $\l_{1}$ contienne $\j_{0}$),
${\underline
\h_{1}}$ par
$\h_{1}=\i_{1}\cap\m_{1}$. En particulier, il faut que
$\i_{p_{0}}$ v\'erifie $\i_{p_{0}}=(
\i_{p_{0}}\cap\m)\oplus (\i_{p_{0}}\cap\a)$.
\end{lem}
\dem On remarque que $I_{2}$,$\dots,I_{p_{0}}$ sont tous  contenus
dans
$ L_{1}\cap L'_{1}$, qui normalise $\n_{1}$.  Comme $I_{1}$
normalise
$\n_{1}$, il en va de m\^{e}me de $x_{1}$. On remarque \'egalement que
${\underline \l_{1}}=Ad\>x_{1}(\l_{1})$, car $\j=Ad\>x_{1}(\j_{0})$.
Par ailleurs
$x_{1}$ est dans $L\cap L'$.  Le lemme est alors une simple
application du transport de structure.\qed

\begin{prop}

Soit $\g$ une alg\`ebre de Lie semi-simple complexe, $\j$ une
\sac de $\g$. On note $\j_{\R}$, l'espace  des \'el\'ements de $\j$,
dont l'image par la repr\'esentation adjointe n'a que des valeurs
propres r\'eelles.
\ste (i) Si $\f\subset \j$ est une \sac fondamentale
d'une forme r\'eelle de $\g$, $\h$,  il existe un  ensemble de
racines
positives de $\Delta(\g, \j)$, $\Delta(\g, \j)^{+}$ une involution
$\C$-lin\'eaire  de
$\j^{*}$, induisant une involution du diagramme de Dynkin
correspondant et telle que :
$$\f=\{X\in \j\vert  ^{t}\tau (X)=-{\overline  X}\}$$
o\`u $X\mapsto  {\overline  X}$ est la conjugaison de $\j$ par
rapport \`a sa forme r\'eelle $\j_{\R}$ et o\`u $^t\tau$ est la
transpos\'ee de
$\tau$.
\ste De plus
$\tau$ est caract\'eris\'ee par cette condition.\ste
(ii) On note $\sigma$ la conjugaison de $\g$, par rapport \`a sa
forme r\'eelle $\h$. Une forme  r\'eelle de $\g$, $\h'$, admet $\f$
pour
\sac, automatiquement fondamentale, si et seulement si elle est
l'espace des points fixes d'une involution antilin\'eaire, $\sigma
'$, de la forme :
$$\sigma'=\sigma\circ Ad\>j, \>\> j\in Z^{1}(\sigma, J)=\{j\in
J\vert j^{\sigma}=j^{-1}\}$$
(iii) Si $j=j_{1}j_{2}j_{2}^{-\sigma}$,on a : $$\sigma\circ Ad\>j=
Ad\>{j_{2}}^{-1}\circ (\sigma\circ Ad\>j_{1})\circ
Ad\>j_{2}$$  (iv) Le quotient $H^{1}(\sigma, J)$ du groupe
$Z^{1}(\sigma, J)$ par son sous-groupe
$B^{1}(\sigma, J)=\{jj^{-\sigma}\vert j\in
J\}$ est fini.
\end{prop}
\dem
On note $\sigma$ la conjugaison de $\g$, par rapport \`a sa
forme r\'eelle $\h$.\ste
Comme $\f$ est une \sac fondamentale de $\h$, on peut choisir un
\'el\'ement r\'egulier de $\g$, $X_{0}$,  dans $\f\cap i\j_{\R}$. On
d\'efinit :
$$\Delta(\g, \j)^{+}:=\{\alpha \in \Delta(\g, \j)\vert i\alpha
(X_{0})>0\} $$
On remarque que la conjugaison complexe induit une involution de
$\Delta(\g, \f)$ car si $\beta\in \Delta(\g, \f)$ et $X\in
\g^{\beta}$, $\sigma(X)$ est \'el\'ement de $\g^{\overline
{\beta}}$.   Si
$\beta\in Hom_{\R}(\f,\C)$,
$\beta$ admet un unique prolongement $\C$-lin\'eaire \`a $\j$. Celui-ci
appartient \`a
$\Delta(\g,
\j)$ et est not\'e ${\tilde \beta}$. Alors l'application :
$$\alpha \mapsto -({\overline {\alpha_{\vert \f}}})^{\tilde {}},
\>\>
\alpha \in \Delta(\g, \j)
$$
est clairement induite  par un endomorphisme $\R$-lin\'eaire du
sous-espace r\'eel de $\j^{*}$ engendr\'e par les \'el\'ements de
$\Delta(\g, \j)$. On note son prolongement $\C$-lin\'eaire \`a
$\j^{*}$,
$\tau$. Alors $\tau $ est involutif. On voit facilement qu'il
pr\'eserve
$\Delta(\g,
\j)^{+}$, et induit donc une involution du diagramme de Dynkin.
\ste Par ailleurs :$$\f=(\f\cap \j_{\R})\oplus (\f\cap i\j_{\R}),
\>\> \j=\f\oplus i\f$$ car $\f$ est une \sac d'une forme r\'eelle de
$\g$, contenue dans
$\j$. Tenant compte de ces d\'ecompositions, on v\'erifie ais\'ement
l'\'egalit\'e :
$$\f=\{X\in \j\vert  ^{t}\tau (X)=-{\overline  X}\}$$
A noter que celle-ci caract\'erise la restriction de $\tau $ \`a $\f$,
donc d\'etermine $\tau$ par $\C$-lin\'earit\'e. D'o\`u (i).
\ste Montrons (ii)). Soit $\h'$ une forme  r\'eelle de $\g$,
admettant
$\f$ pour
\sac. D'abord celle-ci est fondamentale car aucune racine de
$\f$ n'est r\'eelle puisque $\f$ est fondamentale dans $\h$. Notons
$\sigma '$ la conjugaison par rapport \`a $\h'$. Alors $\sigma\sigma
'$ est un automorphisme $\C$-lin\'eaire de l'alg\`ebre de Lie  $\g$,
qui est l'identit\'e sur $\f$, donc sur $\j$, par $\C$-lin\'earit\'e.
Il est donc de la forme $Ad\>j,\>\> j\in J$ (cf. [Bou], Ch. VIII,
Paragraphe 5.2, Proposition 2 ,et Ch. VI).
Donc, on a :
$$\sigma'=\sigma\circ Ad\>j$$
Le fait que $\sigma'$ soit une involution montre que $j$ doit \^{e}tre
\'el\'ement de $ Z^{1}(\sigma, J)$, et $\sigma'$ a la forme voulue.
\ste R\'eciproquement, on remarque que si
$\sigma '$ est de la forme indiqu\'ee, c'est un automorphisme
antilin\'eaire de
$\g$, qui est involutif car $j\in Z^{1}(\sigma, J)$. Son ensemble
de points fixes est une forme r\'eelle de $\g$, $\h'$, contenant $\f$.
Donc $\f$ est une \sac de $\h'$, n\'ecessairement fondamentale,
comme on l'a vu plus haut. Ceci ach\`eve de prouver (ii).\ste
Montrons (iii), en identifiant $J$ \`a $\C^{n}/ \Z^{n}$, la
restriction de
$\sigma$ \`a $J$ \'etant alors d\'efinie par passage au quotient d'un
endomorphisme   antilin\'eaire de $\C^{n}$,  pr\'eservant $\Z^{n}$,
not\'e ${\tilde \sigma}$. Un calcul imm\'ediat  montre que :
$$H^{1}(\sigma, J)\approx (\{X\in \R^{n}\vert \> {\tilde \sigma}
(X)=X, 2X\in
\Z^{n}\}+\Z^{n}) /\Z^{n}$$ Mais le deuxi\`eme membre est contenu
dans $ ((1/2)
\Z^{n})/\Z^{n}$ qui est isomorphe au groupe fini $(\Z/
2\Z)^{n}$. Ceci ach\`eve de prouver (iii) et la
Proposition.\qed
La preuve de la Proposition suivante est imm\'ediate
\begin{prop}
Soit $\g_{1}$ une alg\`ebre de Lie simple complexe, $\j_{1}$ une
\sac de $\g_{1}$. On note $\g:=\g_{1}\times \g_{1}$,
$\j_{0}:=\j_{1}\times \j_{1}$. \ste (i) Un sous-espace $\f$ est
une \sac de l'espace des points fixes d'un automorphisme
involutif $\C$-lin\'eaire de $\g$, permutant les facteurs, si et
seulement si $$\f=\{(X,\tau(X))\vert X\in \g_{1}\}, $$
o\`u $\tau$ est un automorphisme de $\Delta(\g_{1}, \j_{1})$,
l'involution \'etant alors donn\'ee par:
$$\sigma_{{\tilde \tau}}(X,Y)=  ({\tilde \tau}^{-1}(Y),{\tilde
\tau} (X)),
\>\> (X,Y)\in \g_{1}\times \g_{1},$$
o\`u ${\tilde \tau}$ est un automorphisme de $\g_{1}$, dont la
restriction \`a $\j_{1}$ est \'egale \`a $\tau$.\ste
En particulier il n'y a qu'un nombre fini de possibilit\'es pour
$\tau $ et $\f$.\ste   (ii) Si ${\tilde \tau}$ est  un
automorphisme de $\g_{1}$, dont la restriction \`a $\j_{1}$ est
\'egale \`a $\tau$, les automorphismes de $\g_{1}$ poss\'edant la m\^{e}me
propri\'et\'e sont ceux de la forme ${\tilde \tau}\circ Ad\>
j_{1}$, $j_{1}\in J_{1}$.\ste (iii) Les automorphismes involutifs
$\C$-lin\'eaires  de
$\g$ permutant les facteurs, dont $\f$ est une \sac de l'ensemble
de ses points fixes, sont de la forme :
$$Ad\>j^{-1}\circ\sigma_{{\tilde \tau}}\circ Ad\>j$$
o\`u :$$j=(j_{1}, 1), \>\>j_{1}\in J_{1}$$
\end{prop}
\ste  Le Th\'eor\`eme 4 montre que, \`a conjugaison sous
$G$ (ou sous
$P\cap P'$ d'apr\`es la Proposition 2 ), tout triple de Manin
fortement standard de
$\g$, sous
$(\p,\p')$, est obtenu par rel\`evement d'un triple de Manin fortement
standard de $\l\cap\l'$. Il  est clair que l'ensemble des classes
de conjugaison sous $G$ des  rel\`evements d'un
triple de Manin  de
$\l\cap\l'$, ne d\'epend que de la classe de conjugaison sous
$L\cap L'$ de ce triple de Manin, ceci par transport de structure.
C'est pourquoi on peut se limiter \`a appliquer le rel\`evement \`a une
famille ${\mathcal F}$  de repr\'esentants fortement standard des
classes de conjugaison sous $L\cap L'$  de triples de Manin de
$\l\cap\l'$, pour obtenir tous les triples de Manin sous $(\p,\p')$,
\`a conjugaison sous $G$ pr\`es. D'apr\`es la Proposition 4, le rel\`evement
de deux triples de Manin distincts de
${\mathcal F}$ sont non conjugu\'es sous $\g$. La Proposition
suivante permet de donner une condition pour que deux rel\`evements
d'un m\^{e}me \'el\'ement de ${\mathcal F}$ soient conjugu\'es.
\begin{prop}
(i) Si deux triples de Manin de $\g$,  $(B,\i,\i')$,
$(B,{\underline \i},{\underline \i'})$, sous $(\p,\p')$, ayant
 le m\^{e}me ant\'ec\'edent, sont
conjugu\'es sous
$G$, alors il existe $x\in M$, $x'\in M'$ tels que :
$$Ad\>x(\h )={\underline \h}, \>\> Ad\>x'(\h')={\underline
\h'},\>\> o\grave{u} \>\>\h=\i\cap \m,\>\> etc. $$
De plus :$$\i_{\a}={\underline \i}_{\a}\>\> et \>\>
\i'_{\a'}={\underline
\i'_{\a'}}$$
(ii) Si  $(B,\i,\i')$,
$(B,{\underline \i},{\underline \i'})$ sont deux triples de Manin
sous
$(\p,\p')$ pour lequel il existe $x\in M$, $x'\in M'$ tels que :
$$Ad\>x(\i)={\underline \i}, \>\> Ad\>x'(\i')={\underline \i'},$$
ils sont conjugu\'es par un \'el\'ement de $G$, si et seulement si
$x{\tilde I}\cap  x'{\tilde I'}$. Ici ${\tilde I}$  d\'esigne le
groupe ${\tilde H}AN$, o\`u ${\tilde H}=\{m\in M\vert Ad\> m \>\>
pr\acute{e} serve \>\>\h\}$. De plus  ${\tilde I}$ est  le
normalisateur dans $\g$ de $\i$. Le groupe ${\tilde I}'$ est d\'efini
de mani\`ere similaire.
\end{prop}
\dem
Si deux triples de Manin sont comme dans (i), il existe $g\in
P\cap P'$ qui les conjugue, d'apr\`es la Proposition 2 (i). Ecrivons
$g=unn'$, avec
$u\in L\cap L'$,
$n\in N_{L'}$, $n'\in N'_{L}$ et posons ${\tilde x}=un'\in L$,
${\tilde x'}=un\in L'$. Comme $n$ et $n'$ commutent et que $n$
normalise $\i$, on a : $Ad\>{\tilde x}(\i)={\underline \i}$. Comme
${\tilde x}\in L$, $Ad\>{\tilde x}$ pr\'eserve $\m$ et fixe les
\'el\'ements de $\a$. D'o\`u il r\'esulte imm\'ediatement que $Ad\>{\tilde
x}(\h)={\underline
\h}$ et $\i_{\a}={\underline \i}_{\a}$. Ecrivant
${\tilde x}=xa$, avec $x\in M$, $a\in A$, on voit que $x$ v\'erifie
les propri\'et\'es voulues. On proc\`ede de m\^{e}me pour trouver $x'$.
Ceci prouve (i). \ste
Montrons (ii) et soit deux triples de Manin comme dans l'\'enonc\'e.
Ils sont conjugu\'es sous $G$, si et seulement si il existe $y\in G$
tel que :
$$Ad\>y(\i)={\underline \i}, Ad\>y(\i')={\underline \i'}$$
Comme $x\in M$, il normalise $\n$ et fixe $\a$.
Tenant compte des conditions  de l'\'enonc\'e, on voit que  $y$ doit
v\'erifier : $$ Ad\>x^{-1}y(\i)=\i, Ad\>x'^{-1}y(\i')=\i'$$
Pour conclure, il reste \`a prouver que le normalisateur dans $G$ de
$\i$, est bien \'egal \`a ${\tilde I}$. Il est clair que ${\tilde I}$
est bien contenu dans ce normalisateur. R\'eciproquement, soit $y$
un \'el\'ement de ce normalisateur. Il normalise le radical nilpotent
$\n$ de $\i$, donc $y\in P$, d'apr\`es (1.17). On \'ecrit $y=man $ avec
$m\in M$, $a\in A$, $n\in N$. Comme $AN$ est dans le normalisateur
de $\i$, on voit que $m$ normalise $\i$ et donc aussi
$\h=\i\cap\m$. Donc $m\in {\tilde H}$ et $y$ appartient \`a
${\tilde I}$, comme d\'esir\'e. Ceci ach\`eve de prouver la
Proposition.\qed

\ste {\bf R\'ef\'erences }\ste

\noindent[Bor], BOREL A., {\em Linear algebraic groups, Second
Enlarged Edition}, Graduate Text in Math.126, 1991, Springer Verlag,
 New York, Berlin, Heidelberg.

\noindent[Bou], BOURBAKI N., {Groupes et Alg\`ebres de Lie, Chapitre
I, Chapitres IV, V, VI, Chapitres VII, VIII}, Actualit\'es
Scientifiques et Industrielles 1285, 1337, 1364, Hermann, Paris,
1960, 1968, 1975.

\noindent[G], GANTMACHER F.,  {\em
Canonical representation of automorphism of a semisimple Lie
group}, Math Sb., 47, (1939), 101-144.

\noindent[K1], KAROLINSKY E., {\em A classification of Poisson
homogeneous spaces of a compact Poisson Lie group}, Math.
Phys., Anal. and Geom., 3 (1996), 545-563.

\noindent[K2], KAROLINSKY E., {\em A classification of Poisson
homogeneous spaces of a compact Poisson Lie group},
Dokl. Ak. Nauk, 359 (1998), 13-15.

\noindent[K3], KAROLINSKY E.,{\em A classification of Poisson
homogeneous spaces of a reductive complex Poisson Lie group,}
 Preprint, 1999

\noindent[M1], MATSUKI T., {\em The orbits of affine symmetric
spaces under the action of minimal parabolic subgroups}, J. Math.
Soc. Japan, 31 (1979), 331-357.

\noindent[M2], MATSUKI T., {\em Orbits of affine symmetric
spaces under the action of  parabolic subgroups}, Hiroshima J.
Math., 12 (1982), 307-320.

\ste {\em Institut de Math\'ematiques de Luminy, U.P.R. 9016 du
C.N.R.S. \ste
Universit\'e de la M\'editerrann\'ee, 163 Avenue de Luminy, Case
907,\ste
13288, Marseille Cedex 09, France \ste
e-mail : delorme@iml.univ-mrs.fr}

\end{document}